\documentclass[conf]{new-aiaa}

\let\transp\relax

\usepackage{amsmath}
\usepackage{amssymb}
\usepackage{amsthm}
\usepackage{mathtools}
\usepackage{ar}
\usepackage{empheq}
\usepackage{bm} 

\usepackage{color}
\usepackage{xcolor}

\usepackage{ifthen}
\usepackage{pgffor}
\usepackage{etoolbox}
\usepackage{parselines}
\usepackage{environ}
\usepackage{listofitems}
\usepackage{forloop}
\usepackage{calc}

\usepackage{graphicx}
\usepackage{subcaption}

\usepackage{tabularx}
\usepackage{longtable}
\usepackage{multicol}
\usepackage{diagbox} 
\usepackage{varwidth} 
\usepackage{makecell} 
\usepackage{arydshln}

\usepackage{tikz}
\usepackage[tikz]{mdframed}
\usepackage{enumitem}
\usepackage{algorithm, algorithmicx}
\usepackage{algpseudocode}
\usepackage[many]{tcolorbox}

\usepackage{float}
\usepackage{afterpage}
\usepackage{pdflscape} 

\usepackage{hyperref}

\usepackage{xspace}

\usepackage{varwidth}
\usepackage{siunitx}
\usepackage{nomencl}
\usepackage[symbol*]{footmisc}

\makeatletter
\newcommand\footnoteref[1]{\protected@xdef\@thefnmark{\ref{#1}}\@footnotemark}
\makeatother

\newcommand{\definedas}{\coloneqq}
\newcommand{\fr}[1]{\mathcal{F}_{\mathcal{#1}}}
\newcommand{\integer}{\mathbb Z}
\newcommand{\quaternion}{\mathbb Q}
\newcommand{\ones}[1]{\bm{1}_#1}
\newcommand{\blkdiag}[1]{\mathrm{blkdiag}(#1)}

\usepackage{rotating}
\usepackage{tikz}
\usepackage{pgfplots}
\usepackage{ifthen}
\pgfplotsset{compat=1.12}
\usetikzlibrary{calc}
\usetikzlibrary{patterns}
\usetikzlibrary{decorations.pathmorphing,decorations.markings}
\usetikzlibrary{math}
\usetikzlibrary{scopes}
\usetikzlibrary{fadings}
\usetikzlibrary{arrows,bending}
\usetikzlibrary{arrows.meta}
\tikzset{>=latex}

\pgfdeclarelayer{bg}    
\pgfsetlayers{bg,main}  

\definecolor{beige}{RGB}{245,245,220}
\definecolor{darkred}{rgb}{0.90,0.00,0.00}
\definecolor{darkgreen}{rgb}{0.00,0.45,0.00}

\definecolor{ucol}{RGB}{255,0,0}
\definecolor{gcol}{RGB}{0,120,0}
\definecolor{scol}{RGB}{63, 226, 45}
\definecolor{vcol}{RGB}{0,0,0}
\definecolor{acol}{RGB}{0, 128, 255}
\definecolor{dcol}{RGB}{204,102,0}
\definecolor{lcol}{RGB}{204,102,0}
\definecolor{bcol}{RGB}{0,0,0}
\definecolor{ocol}{RGB}{167,167,167}
\definecolor{propcolor}{rgb}{1.0,0.0,0.0}
\definecolor{solvecolor}{rgb}{0.0,0.55,0.0}
\definecolor{scalecolor}{RGB}{75, 47, 132}
\definecolor{fillcolor123}{rgb}{1,1,1}
\definecolor{chromeyellow}{rgb}{1.0, 0.65, 0.0}
\definecolor{engine}{rgb}{0.28, 0.24, 0.2}
\definecolor{bodygray}{rgb}{0.28, 0.24, 0.2}
\definecolor{hatch}{rgb}{0.59, 0.0, 0.09}
\colorlet{bodycol}{white!60!black}
\colorlet{darkgrey}{white!40!black}
\definecolor{redengine}{RGB}{152,51,0}
\colorlet{rcscol}{bodycol!75!black}
\colorlet{cmcol}{white!50!black}
\definecolor{coast}{rgb}{0.13, 0.55, 0.13}
\definecolor{burn}{rgb}{0.89, 0.0, 0.13}
\definecolor{panecol}{rgb}{0.85, 0.65, 0.13}
\definecolor{vennc1}{HTML}{3F3DBA}
\definecolor{vennc2}{HTML}{3DBA3F}
\definecolor{vennc3}{HTML}{B8BA3D}
\definecolor{vennc4}{HTML}{BA3DB8}
\definecolor{vennfill}{HTML}{A8B6BF}
\colorlet{spcol}{blue!35!white!95!black}
\definecolor{Xcol}{RGB}{255,0,0}
\definecolor{Ycol}{RGB}{128,128,0}
\definecolor{Zcol}{RGB}{0,0,255}

\tikzstyle{block_blank} = [draw, thick, rectangle, minimum height=3em, minimum width=3em]
\tikzstyle{block} = [draw, ultra thick, fill=blue!20, rectangle, 
    rounded corners=6pt, minimum height=3em, minimum width=6em]
\tikzstyle{mexblock} = [draw, ultra thick, magenta, fill=blue!20, rectangle, 
    rounded corners=6pt, minimum height=3em, minimum width=6em, text=black]
\tikzstyle{sum} = [draw, fill=blue!20, circle, node distance=2cm]
\tikzstyle{sum_blank} = [draw, circle, node distance=2cm]
\tikzstyle{input} = [coordinate]
\tikzstyle{output} = [coordinate]
\tikzstyle{pinstyle} = [pin edge={to-,thin,black}]

\newcommand{\threeaxes}[8]{
	\tikzmath{
    	\Lxr=#3;\Lxl=#3;\Lyt=#4;\Lyb=#4;\Zt= #5;\Zb= #6;\Xang=#7;\Yang=#8;
    	\Xxr= cos(\Xang)*\Lxr; \Xyr=-sin(\Xang)*\Lxr;
    	\Xxl=-cos(\Xang)*\Lxl; \Xyl= sin(\Xang)*\Lxl;
        \Yxt= sin(\Yang)*\Lyt; \Yyt= cos(\Yang)*\Lyt;
        \Yxb=-sin(\Yang)*\Lyb; \Yyb=-cos(\Yang)*\Lyb;
        \zzz=0;
    }
    \begin{scope}[shift={(#1,#2)},rotate=0]
        \ifx\Zt\zzz\else  \draw[black,->] (0,0) -- +(0, \Zt);    \fi
        \ifx\Zb\zzz\else  \draw[black,->] (0,0) -- +(0,-\Zb);    \fi
        \ifx\Lxr\zzz\else \draw[black,->] (0,0) -- +(\Xxr,\Xyr); \fi
        \ifx\Lxl\zzz\else \draw[black,->] (0,0) -- +(\Xxl,\Xyl); \fi
        \ifx\Lyt\zzz\else \draw[black,->] (0,0) -- +(\Yxt,\Yyt); \fi
        \ifx\Lyb\zzz\else \draw[black,->] (0,0) -- +(\Yxb,\Yyb); \fi
	\end{scope}
}
\newcommand{\threeaxeslabelx}[9]{
	\tikzmath{
    	\Lxr=#3;\Lxl=#3;\Lyt=#4;\Lyb=#4;\Zt= #5;\Zb= #6;\Xang=#7;\Yang=#8;
    	\Xxr= cos(\Xang)*\Lxr; \Xyr=-sin(\Xang)*\Lxr;
    	\Xxl=-cos(\Xang)*\Lxl; \Xyl= sin(\Xang)*\Lxl;
        \Yxt= sin(\Yang)*\Lyt; \Yyt= cos(\Yang)*\Lyt;
        \Yxb=-sin(\Yang)*\Lyb; \Yyb=-cos(\Yang)*\Lyb;
        \zzz=0;
    }
    \begin{scope}[shift={(#1,#2)},rotate=0]
        \ifx\Lxr\zzz\else \draw (\Xxr,\Xyr) node[anchor=south west] {#9}; \fi
	\end{scope}
}
\newcommand{\threeaxeslabely}[9]{
	\tikzmath{
    	\Lxr=#3;\Lxl=#3;\Lyt=#4;\Lyb=#4;\Zt= #5;\Zb= #6;\Xang=#7;\Yang=#8;
    	\Xxr= cos(\Xang)*\Lxr; \Xyr=-sin(\Xang)*\Lxr;
    	\Xxl=-cos(\Xang)*\Lxl; \Xyl= sin(\Xang)*\Lxl;
        \Yxt= sin(\Yang)*\Lyt; \Yyt= cos(\Yang)*\Lyt;
        \Yxb=-sin(\Yang)*\Lyb; \Yyb=-cos(\Yang)*\Lyb;
        \zzz=0;
    }
    \begin{scope}[shift={(#1,#2)},rotate=0]
        \ifx\Lyt\zzz\else \draw (\Yxt,\Yyt) node[anchor=south west] {#9}; \fi
	\end{scope}
}
\newcommand{\threeaxeslabelz}[9]{
	\tikzmath{
    	\Lxr=#3;\Lxl=#3;\Lyt=#4;\Lyb=#4;\Zt= #5;\Zb= #6;\Xang=#7;\Yang=#8;
    	\Xxr= cos(\Xang)*\Lxr; \Xyr=-sin(\Xang)*\Lxr;
    	\Xxl=-cos(\Xang)*\Lxl; \Xyl= sin(\Xang)*\Lxl;
        \Yxt= sin(\Yang)*\Lyt; \Yyt= cos(\Yang)*\Lyt;
        \Yxb=-sin(\Yang)*\Lyb; \Yyb=-cos(\Yang)*\Lyb;
        \zzz=0;
    }
    \begin{scope}[shift={(#1,#2)},rotate=0]
        \ifx\Zt\zzz\else  \draw (0,\Zt)     node[anchor=west] {$\;$#9}; \fi
	\end{scope}
}

\newcommand{\isoaxesNL}[6]{
    \tikzmath{
        \rot=#3;
        \len=#4;
        \ddd=#5;
        \pX=  0; \Xx=cos(\pX)*\len; \Xy=sin(\pX))*\len;
        \pY=120; \Yx=cos(\pY)*\len; \Yy=sin(\pY))*\len;
        \pZ=240; \Zx=cos(\pZ)*\len; \Zy=sin(\pZ))*\len;
        \Xxx=cos(\ddd)*\Xx-sin(\ddd)*\Xy; \Xxy=sin(\ddd)*\Xx+cos(\ddd)*\Xy;
        \Yxx=cos(\ddd)*\Yx-sin(\ddd)*\Yy; \Yxy=sin(\ddd)*\Yx+cos(\ddd)*\Yy;
        \Zxx=cos(\ddd)*\Zx-sin(\ddd)*\Zy; \Zxy=sin(\ddd)*\Zx+cos(\ddd)*\Zy;
    }
    \begin{scope}[shift={(#1,#2)},rotate=\rot]
		\filldraw[black] (0,0) circle (2pt);
        \draw[Xcol,thick,->] (0,0) -- +(\Xx,\Xy);
        \draw[Ycol,thick,->] (0,0) -- +(\Yx,\Yy);
        \draw[Zcol,thick,->] (0,0) -- +(\Zx,\Zy);
        \draw (0.1,0.6)   node[rotate=0,anchor=center] {#6};
    \end{scope}
}

\newcommand{\pane}[7]{
	\tikzmath{
    	\rot=#3;
    	\width=#4;
        \height=#5;
        \corner=#6;
    	\px1= 0.5*\width-\corner; \py1=-0.5*\height;
        \px2= 0.5*\width;         \py2=-0.5*\height+\corner;
        \px3= 0.5*\width;         \py3= 0.5*\height-\corner;
        \px4= 0.5*\width-\corner; \py4= 0.5*\height;
        \px5=-0.5*\width+\corner; \py5= 0.5*\height;
        \px6=-0.5*\width;         \py6= 0.5*\height-\corner;
        \px7=-0.5*\width;         \py7=-0.5*\height+\corner;
        \px8=-0.5*\width+\corner; \py8=-0.5*\height;
    }
	\begin{scope}[shift={(#1,#2)},rotate=\rot]
    	\filldraw[{#7}]
        (\px1,\py1) to[out=    0,in=  -90] (\px2,\py2) --
        (\px3,\py3) to[out=   90,in=    0] (\px4,\py4) --
        (\px5,\py5) to[out= -180,in=   90] (\px6,\py6) --
        (\px7,\py7) to[out=  -90,in= -180] (\px8,\py8) -- cycle;
    \end{scope}
}

\newcommand{\coneback}[7]{
	\tikzmath{\rot=#3;
              \length=#4; \radius=\length*tan(0.5*#5); \depth=#6;
              \sx =  cos(\rot)*#1 + sin(\rot)*#2;
              \sy = -sin(\rot)*#1 + cos(\rot)*#2;
    }
    \begin{scope}[shift={(\sx,\sy)},transform canvas={rotate=\rot}]
	    \draw[{#7}] (\radius,-\length) arc(360:180: {\radius} and {-\depth});
    \end{scope}
}

\newcommand{\cone}[7]{
	\tikzmath{\rot=#3;
              \length=#4; \radius=\length*tan(0.5*#5); \depth=#6;
              \sx =  cos(\rot)*#1 + sin(\rot)*#2;
              \sy = -sin(\rot)*#1 + cos(\rot)*#2;
	}
    \begin{scope}[shift={(\sx,\sy)},transform canvas={rotate=\rot}]
    	\fill[{#7}] (0,0) -- (\radius,-\length) arc(360:180: {\radius} and {\depth}) -- cycle;
    	\draw[color=black!100] (0,0) -- (\radius,-\length) arc(360:180: {\radius} and {\depth}) -- cycle;
    \end{scope}
}

\newcommand{\apollodocked}[6]{
\tikzmath{
    \cx=#1; \cy=#2; \sc=#3; \aa=#4; \bb=#5; \throttle=#6;
    \dx = 3.5; \dy = 1.5; \ddx=3; \ddy=0.8;
    \fx=0.4; \ye=-0.75; \we=0.25; \wwe=1.5; \he=1.5;
    \tx=0.1;
    \dddx=0.5; \dddy=0.45;
    \cc=0.35; \comrad=0.4; \comradi=0.3*\comrad; \comradd=\comrad-\comradi;
    \raxi=sin(45)*\comradi; \rayi=cos(45)*\comradi;
    \raxd=sin(45)*\comradd; \rayd=cos(45)*\comradd;
    \sx=(cos(\aa)*\cx + sin(\aa)*\cy)/\sc; 
    \sy=(-sin(\aa)*\cx + cos(\aa)*\cy)/\sc; 
    }

\begin{scope}[shift={(\sx,\sy)},transform canvas={rotate=\aa,scale=\sc}]
    \coordinate (cg) at (0,0);
    \coordinate (A) at (0,-0.5);
    \coordinate (A2) at (0,-1);
    \coordinate (dim) at (\dx,\dy);
    \coordinate (dim2) at (\ddx,\ddy);
    \coordinate (UR) at ($(A)+0.5*(dim)$);
    \coordinate (LL) at ($(A)-0.5*(dim)$);
    \coordinate (LR) at ($(A)+0.5*(\dx,-\dy)$);
    \coordinate (UL) at ($(A)+0.5*(-\dx,\dy)$);
    \coordinate (UR2) at ($(A2)+0.5*(dim2)$);
    \coordinate (LL2) at ($(A2)-0.5*(dim2)$);
    \coordinate (ULbox) at ($(cg)+(\dddx,\dddy)$);
    \coordinate (LRhatch) at ($(ULbox)+(-0.15,0.1)$);
    \coordinate (UMbody) at (0,3.5);
    \coordinate (LMbody) at (0,0.75);
    \coordinate (COM1) at ($(UMbody)+(0,0.4)$);
    \coordinate (COM2) at ($(COM1)+(-1.8,0)$);
    \coordinate (spire) at (1,3.5);
    \coordinate (RCSL) at (-1.6,1.6);
    \begin{scope}[shift={(0,\ye)},rotate=\bb]
        \fill[fill=red!100,shading=axis,shading angle=90,left color=red!75!orange,right color=orange!80!yellow]
                (-0.5*\wwe+\tx,-\he) -- (0.5*\wwe-\tx,-\he) -- (0,-\he-2*\throttle) -- cycle;
        \draw[fill=black!40,thick,shading=ball,left color=black!40, right color=black!60] (-0.5*\we,0) -- (0.5*\we,0) to[out=-40,in=90] (0.5*\wwe,-\he) -- (-0.5*\wwe,-\he) to[out=90,in=220] cycle;
        \draw [color=black!75,thin] (0.5*\wwe-0.1,-\he+0.2) -- (-0.5*\wwe+0.1,-\he+0.2);
        \draw [color=black!70,thin] (0.5*\wwe-0.15,-\he+0.4) -- (-0.5*\wwe+0.15,-\he+0.4);
        \draw [color=black!60,thin] (0.5*\wwe-0.2,-\he+0.6) -- (-0.5*\wwe+0.2,-\he+0.6);
        \draw [color=black!55,thin] (0.5*\wwe-0.25,-\he+0.75) -- (-0.5*\wwe+0.25,-\he+0.75);
    \end{scope}
    \draw [fill=bodygray,thick] (LL2) rectangle (UR2);
    \draw [color=chromeyellow!50!black, fill=chromeyellow,thick,shading=axis,shading angle=90, left color=chromeyellow,right color=chromeyellow!75!black] (LL) rectangle (UR);
    \coordinate (CCR) at (2.2,-1.2);
    \draw [fill=chromeyellow!75!black] (UR) -- (3,0) -- (3,-0.1) -- ($(UR)-(0,0.1)$) -- cycle;
    \draw [fill=chromeyellow!75!black] ($(LR)+(0,0.2)$) -- (2.98,0) -- (3,-0.1) -- ($(LR)+(0,0.1)$) -- cycle;
    \draw [fill=chromeyellow!75!black] ($(LR)+(0,0.15)$) -- (CCR) -- ($(CCR)+(0,-0.05)$) -- ($(LR)+(0,0.1)$) -- cycle;
    \draw [fill=chromeyellow!75!black] (3,0) -- (CCR) -- (1.27,-2.6) -- (1.36,-2.6) -- (3,-0.1) -- cycle;
    \begin{scope}[shift={(1.33,-2.55)},rotate=-45]
    \draw [fill=chromeyellow!50!black] (-\fx,0) -- (\fx,0) to[bend left] cycle;
    \end{scope}
    \coordinate (CCL) at (-2.2,-1.2);
    \draw [fill=chromeyellow] (UL) -- (-3,0) -- (-3,-0.1) -- ($(UL)-(0,0.1)$) -- cycle;
    \draw [fill=chromeyellow] ($(LL)+(0,0.2)$) -- (-2.98,0) -- (-3,-0.1) -- ($(LL)+(0,0.1)$) -- cycle;
    \draw [fill=chromeyellow] ($(LL)+(0,0.15)$) -- (CCL) -- ($(CCL)+(0,-0.05)$) -- ($(LL)+(0,0.1)$) -- cycle;
    \draw [fill=chromeyellow] (-3,0) -- (CCL) -- (-1.27,-2.6) -- (-1.36,-2.6) -- (-3,-0.1) -- cycle;
    \begin{scope}[shift={(-1.33,-2.55)},rotate=45]
    \draw [fill=chromeyellow!50!black] (-\fx,0) -- (\fx,0) to[bend left] cycle;
    \end{scope}
    \draw [fill=chromeyellow] (0.1,0.25) -- (-0.1,0.25) -- (-0.1,0) -- (0.1,0) -- cycle;
    \draw [fill=chromeyellow] (-0.1,0) -- (0.1,0) -- (0.1,-2.5) -- (-0.1,-2.5) -- cycle;
    \draw [fill=chromeyellow!50!black] (-\fx,-2.5) -- (\fx,-2.5) to[bend left] cycle;
    \draw [fill=bodygray!80,thick,shading=axis,left color=bodygray!60,right color=bodygray!90] (UMbody) -- ++(1.2,0) -- ++(\cc,-\cc) -- ++(0,-1.5) -- ++(0.5,0) -- ++(0.3,-0.2) -- ++(0,-0.5) -- ++(-0.25*\cc,-0.5*\cc) -- ++(-\cc,-0.75*\cc) -- ++(-1,0) -- (LMbody) -- ++(-0.75,-0.35) -- ++(-0.6,0) -- ++(-0.55,0.5) -- ++(\cc,0) -- ++(0,2.25) -- ++(\cc,\cc) -- cycle;
    \draw [fill=bodygray!80,thick] (-1.05,0.25) -- ++(0.1,0) -- ++(0,0.15) -- ++(-0.1,0) -- cycle;
    \draw [fill=bodygray!80,thick] (1.2,0.25) -- ++(0.1,0) -- ++(0,0.25) -- ++(-0.1,0) -- cycle;
    \draw [fill=bodygray!70,thick,shading=ball,left color=bodygray!80,right color=bodygray!90] (0,1.8) circle (1.4); 
    \draw [fill=bodygray,thick] (ULbox) -- ++(-2*\dddx,0) -- ++(0,-0.4) -- ++(1,0) -- cycle;
    \draw [fill=bodygray!80,thick] (ULbox) -- ++(-2*\dddx,0) -- ++(0.2,3) -- ++(0.6,0) -- cycle;
    \draw [fill=hatch!80!black] (LRhatch) -- ++(-2*\dddx+0.3,0) -- ++(-0.09,0.1) -- ++(0.04,0.6) -- ++(0.11,0.11) -- ++(0.6,0) -- ++(0.09,-0.1) -- ++(0.04,-0.6) -- cycle;
    \draw [color=black!80!hatch] ($(LRhatch)+(0,0.1)$) -- ++(-0.7,0);
\end{scope}
}

\newcommand{\apolloCSM}[5]{
\tikzmath{
    \cx=#1/#5; \cy=#2/#5; \aa=#3; \gim=#4; 
    \Lx=3; \Ly=4; \dx=0.04; \dy=0.15; \ddx=\dx; \ddy=0.05;
    \engx=0; \engy=0; \we=0.25; \wwe=1.6; \he=1.75; \hee=0.5*\he;
    \wpx=-1; \wpy=0.8; \lwp=2.48;
    \twpx=-0.85; \twpy=3.8; \twpw=1; \twph=0.4; \twpxx=\twpx+\twpw+0.25;
    \tbpx=\twpx-2*0.25; \tbpy=\twpy; \tbpw=0.3; \tbph=\twph;
    \rcx=0; \rcy=2; \rcdx=0.25; \rcdy=0.45; \rch=0.25; 
    \rcxx=-0.5*\Lx; \rcyy=\rcy; \rcdxx=0.1; 
    \rcxxx=0.5*\Lx; \rcyyy=\rcy; 
    \ax=1; \ady=0.05; \ay=0.5*\ady; \adx=1.25; 
    \cmx=0; \cmy=\Ly+\ddy; \cmh=0.45*\Ly; \cmw=0.3*\Lx;
}

\begin{scope}[shift={(\cx,\cy)},rotate=\aa,transform canvas={scale=#5}]
    \begin{scope}[shift={(\engx,\engy)},rotate=\gim]
        \coordinate (EUR) at (0.5*\we,0);
        \coordinate (EUL) at (-0.5*\we,0);
        \draw[fill=black!40,thick,shading=ball,left color=black!40, right color=black!60] (EUL) -- (EUR) to[out=-40,in=95] (0.5*\wwe,-\he) -- (-0.5*\wwe,-\he) to[out=85,in=220] cycle;
        \draw[fill=redengine,shading=ball,left color=redengine!90!white!80!black, right color=redengine!55!black] (EUL) -- (EUR) to[out=-40,in=105] (0.42*\wwe,-\hee) -- (-0.42*\wwe,-\hee) to[out=75,in=220] cycle;
    \end{scope}
    \coordinate (LL) at (-0.5*\Lx,0);
    \coordinate (LR) at (+0.5*\Lx,0);
    \coordinate (UR) at (+0.5*\Lx,\Ly);
    \coordinate (UL) at (-0.5*\Lx,\Ly);
    \draw[fill=darkgrey,rounded corners=3pt] ($(LL)+(\dx,\dy)$) -- ($(LR)+(-\dx,\dy)$) -- ($(LR)+(-\dx,-\dy)$) -- ($(LL)+(\dx,-\dy)$) -- cycle;
    \draw[fill=darkgrey,rounded corners=2pt] ($(UL)+(\ddx,\ddy)$) -- ($(UR)+(-\ddx,\ddy)$) -- ($(UR)+(-\ddx,-\ddy)$) -- ($(UL)+(\ddx,-\ddy)$) -- cycle;
    \draw[thick, rounded corners=2pt,fill=bodycol,shading=axis, left color=white!90!black, right color=bodycol] (LL) -- (LR) -- (UR) -- (UL) -- cycle; %
    \draw[opacity=0.8,shading=axis,right color=white!85!black, left color=white] (\wpx,\wpy) -- (\wpx+\lwp,\wpy) -- (\wpx+\lwp,\dy) -- (\wpx,\dy) -- cycle;
    \foreach \i in {2,3,4,5,6}{
        \draw[opacity=0.4] (\wpx+\ddx,\i*0.12) -- (\wpx+\lwp-\ddx,\i*0.12);
    }
    \draw[opacity=0.8,shading=axis,right color=white!90!black, left color=white] (\twpx,\twpy) -- ++(\twpw,0) -- ++(0,-\twph) -- ++(-\twpw,0) -- cycle;
    \foreach \i in {1,2,3}{
        \draw[opacity=0.4] (\twpx+\dx,\twpy-\i*0.1) -- ++(\twpw-2*\dx,0);
    }
    \draw[opacity=0.7,shading=axis,right color=white!80!black, left color=white] (\twpxx,\twpy) -- ++(\twpw,0) -- ++(0,-\twph) -- ++(-\twpw,0) -- cycle;
    \foreach \i in {1,2,3}{
        \draw[opacity=0.4] (\twpxx+\dx,\twpy-\i*0.1) -- ++(\twpw-2*\dx,0);
    }
    \draw[opacity=0.8,shading=axis,right color=bodycol!70!black,left color=bodycol!90!black] (\tbpx,\tbpy) -- ++(\tbpw,0) -- ++(0,-\tbph) -- ++(-\tbpw,0) -- cycle;
    \draw[opacity=0.8] (\tbpx+0.02,\tbpy-0.02) -- ++(0,-\tbph+0.04) -- ++(\tbpw-0.04,0) -- ++(0,0.9*\tbph) -- cycle;
    \draw[opacity=0.8] (\tbpx+0.1,\tbpy-0.1) circle (0.75pt);
    \draw[opacity=0.8] (\tbpx+0.2,\tbpy-0.1) circle (0.75pt);
    \draw[opacity=0.8] (\tbpx+0.15,\tbpy-0.2) circle (0.75pt);
    \draw[opacity=0.8] (\tbpx+0.1,\tbpy-0.3) circle (0.75pt);
    \draw[opacity=0.8] (\tbpx+0.2,\tbpy-0.3) circle (0.75pt);

    \begin{scope}[shift={(\rcx,\rcy)}]
        \draw[] (0.5*\rcdx,0.5*\rcdy) -- ++(0,-\rcdy) -- ++(-\rcdx,0) -- ++(0,\rcdy) -- cycle;
        \draw[fill=bodycol!60!black] (0.5*\rcdx,0.5*\rcdy) -- (0.35*\rcdx,0.35*\rcdy) -- ++(0,-0.7*\rcdy) -- (0.5*\rcdx,-0.5*\rcdy);
        \draw[fill=bodycol!90!black] (-0.5*\rcdx,0.5*\rcdy) -- (-0.35*\rcdx,0.35*\rcdy) -- (0.35*\rcdx,0.35*\rcdy) -- (0.5*\rcdx,0.5*\rcdy);
        \draw[fill=bodycol!90!black] (-0.5*\rcdx,-0.5*\rcdy) -- (-0.35*\rcdx,-0.35*\rcdy) -- (0.35*\rcdx,-0.35*\rcdy) -- (0.5*\rcdx,-0.5*\rcdy);
        \draw[fill=bodycol] (-0.5*\rcdx,-0.5*\rcdy) -- (-0.35*\rcdx,-0.35*\rcdy) -- (-0.35*\rcdx,0.35*\rcdy) -- (-0.5*\rcdx,0.5*\rcdy);
        \draw[fill=rcscol,shading=axis,right color=rcscol!75!black] (0.1*\rcdx,-0.45*\rcdy) -- (-0.1*\rcdx,-0.45*\rcdy) to[out=240,in=85] ++(-0.3*\rcdx,-\rch) -- ++(0.8*\rcdx,0) to[out=95,in=-60] cycle;
        \draw[fill=rcscol,shading=axis,right color=rcscol!75!black] (0.1*\rcdx,0.45*\rcdy) -- (-0.1*\rcdx,0.45*\rcdy) to[out=120,in=275] ++(-0.3*\rcdx,\rch) -- ++(0.8*\rcdx,0) to[out=265,in=60] cycle;
        \draw[fill=rcscol,shading=axis,left color=rcscol!80!black, right color=rcscol!60!black] (0.45*\rcdx,0.25*\rcdx) -- (0.45*\rcdx,0.05*\rcdx) to[out=-30,in=175] ++(\rch,-0.3*\rcdx) -- ++(0,0.8*\rcdx) to[out=185,in=30] cycle;
        \draw[fill=rcscol] (-0.45*\rcdx,-0.05*\rcdx) -- (-0.45*\rcdx,-0.25*\rcdx) to[out=210,in=5] ++(-\rch,-0.3*\rcdx) -- ++(0,0.8*\rcdx) to[out=-5,in=150] cycle;
    \end{scope}

    \begin{scope}[shift={(\rcxx,\rcyy)}]
        \draw[fill=rcscol,shading=axis,left color=rcscol!70!white,right color=rcscol!80!black] (0.5*\rcdxx,0.5*\rcdy) -- ++(-\rcdxx,0) -- ++(0,-\rcdy) -- ++(\rcdxx,0) -- cycle;
        \begin{scope}[shift={(0.04,0)},rotate=-10]
        \draw[fill=rcscol,shading=axis,left color=rcscol!70!white,right color=rcscol!80!black] (0.3*\rcdxx,-0.45*\rcdy) -- (-0.3*\rcdxx,-0.45*\rcdy) to[out=240,in=85] ++(-0.5*\rcdxx,-\rch) -- ++(1.6*\rcdxx,0) to[out=95,in=-60] cycle;
        \end{scope}
        \begin{scope}[shift={(0.04,0)},rotate=10]
        \draw[fill=rcscol,shading=axis,left color=rcscol!70!white,right color=rcscol!80!black] (0.3*\rcdxx,0.45*\rcdy) -- (-0.3*\rcdxx,0.45*\rcdy) to[out=120,in=275] ++(-0.5*\rcdxx,\rch) -- ++(1.6*\rcdxx,0) to[out=265,in=60] cycle;
        \end{scope}
        \draw[fill=rcscol,shading=ball,outer color=rcscol, inner color=black] (0,0.05) circle (0.8*\rcdxx);
    \end{scope}
    \begin{scope}[shift={(\rcxxx,\rcyyy)}]
        \draw[fill=rcscol,shading=axis,right color=rcscol!80!black] (0.5*\rcdxx,0.5*\rcdy) -- ++(-\rcdxx,0) -- ++(0,-\rcdy) -- ++(\rcdxx,0) -- cycle;
        \begin{scope}[shift={(-0.04,0)},rotate=10]
        \draw[fill=rcscol,shading=axis,right color=rcscol!80!black] (0.3*\rcdxx,-0.45*\rcdy) -- (-0.3*\rcdxx,-0.45*\rcdy) to[out=240,in=85] ++(-0.5*\rcdxx,-\rch) -- ++(1.6*\rcdxx,0) to[out=95,in=-60] cycle;
        \end{scope}
        \begin{scope}[shift={(-0.04,0)},rotate=-10]
        \draw[fill=rcscol,shading=axis,right color=rcscol!80!black] (0.3*\rcdxx,0.45*\rcdy) -- (-0.3*\rcdxx,0.45*\rcdy) to[out=120,in=275] ++(-0.5*\rcdxx,\rch) -- ++(1.6*\rcdxx,0) to[out=265,in=60] cycle;
        \end{scope}
        \draw[fill=rcscol,shading=ball,outer color=rcscol, inner color=black] (0,-0.05) circle (0.8*\rcdxx);
    \end{scope}
    \begin{scope}[shift={(\ax,\ay)}]
        \draw[fill=darkgrey] (0,0) -- ++(\adx,0) -- ++(0,-\ady) -- ++(-\adx,0) -- cycle;
        \draw[fill=darkgrey] (0.55*\adx,\ady) -- ++(0.1*\adx,0) -- ++(0,-\ady) -- ++(-0.1*\adx,0) -- cycle;
        \draw[fill=darkgrey!80!black] (0.45*\adx,2.5*\ady) to[out=-60,in=240] ++(0.3*\adx,0) -- cycle;
        \draw[fill=darkgrey] (0.55*\adx,2.5*\ady) -- +(0.05*\adx,0.5*\ady) -- +(0.1*\adx,0) -- cycle;
        \draw[fill=darkgrey] (0.9*\adx,\ady) -- ++(0.1*\adx,0) -- ++(0,-\ady) -- ++(-0.1*\adx,0) -- cycle;
        \draw[fill=darkgrey!80!black] (0.8*\adx,2.5*\ady) to[out=-60,in=240] ++(0.3*\adx,0) -- cycle;
        \draw[fill=darkgrey] (0.9*\adx,2.5*\ady) -- +(0.05*\adx,0.5*\ady) -- +(0.1*\adx,0) -- cycle;
    \end{scope}
    \apolloCM{\cmx}{\cmy}{\Lx}{\Ly}
\end{scope}
}

\newcommand{\apolloCM}[4]{
\tikzmath{
    \cmx=#1; \cmy=#2;
    \cmw=0.3*#3; \cmh=0.45*#4;
    \cmtx=-0.9; \cmty=0.15;
    \spx=0; \spy=\cmh; \spw=0.1; \sph=0.4;
}
\begin{scope}[shift={(\cmx,\cmy)}]
        \coordinate (LMUR) at (0.5*\cmw,\cmh);
        \coordinate (LMUL) at (-0.5*\cmw,\cmh);
        \draw[rounded corners=1pt,fill=cmcol,path fading=west] (-0.5*\Lx+\ddx,0) -- (0.5*\Lx-\ddx,0) -- (LMUR) -- (LMUL) -- cycle;
        \draw[rounded corners=1pt] (-0.5*\Lx+\ddx,0) -- (0.5*\Lx-\ddx,0) -- (LMUR) -- (LMUL) -- cycle;
        \draw[opacity=0.1] (-0.85*\cmw,0.65*\cmh) -- (0.85*\cmw,0.65*\cmh);
        \draw[opacity=0.1] (-1.35*\cmw,0.2*\cmh) -- (1.35*\cmw,0.2*\cmh);
        \draw[fill=cmcol!80!white,rounded corners=2pt] (-0.225,0.6*\cmh) -- ++(0.45,0) -- ++(0.175,-0.6) -- ++(-0.8,0) -- cycle;
        \draw[fill=blue!40!black,rounded corners=1pt] (-0.1,0.525*\cmh) -- ++(0.2,0) -- ++(0.05,-0.3) -- ++(-0.3,0) -- cycle;
        \begin{scope}[shift={(\cmtx,\cmty)},rotate=45]
            \draw[fill=black,shading=ball,inner color=black,outer color=white!40!black] (0,0) circle (0.1 and 0.04);
        \end{scope}
        \begin{scope}[shift={(\cmtx-0.2,\cmty)},rotate=45]
            \draw[fill=black,shading=ball,inner color=black,outer color=white!40!black] (0,0) circle (0.04 and 0.1);
        \end{scope}
        \begin{scope}[shift={(-\cmtx,\cmty)},rotate=45]
            \draw[fill=black,shading=ball,inner color=black,outer color=white!40!black] (0,0) circle (0.04 and 0.1);
        \end{scope}
        \begin{scope}[shift={(-\cmtx+0.2,\cmty)},rotate=45]
            \draw[fill=black,shading=ball,inner color=black,outer color=white!40!black] (0,0) circle (0.1 and 0.04);
        \end{scope}
        \begin{scope}[shift={(\spx,\spy)}]
            \draw[fill=cmcol!80!white] (-0.5*\spw,0) -- ++(\spw,0) -- ++(0,\sph) -- ++(0.5*\spw,0) -- ++(-\spw,0.3*\sph) -- ++(-\spw,-0.3*\sph) -- ++(0.5*\spw,0) -- cycle;
            \draw[fill=cmcol!80!white,rounded corners=0.5pt] (0.5*\spw,0) -- (0.5*\spw+0.45*\sph,0.45*\sph) -- (0.5*\spw,0.9*\sph) -- ++(0,-0.5*\spw) -- (0.45*\sph,0.45*\sph) -- (0.5*\spw,0.5*\spw) -- cycle;
            \draw[fill=cmcol!80!white,rounded corners=0.5pt] (-0.5*\spw,0) -- (-0.5*\spw-0.45*\sph,0.45*\sph) -- (-0.5*\spw,0.9*\sph) -- ++(0,-0.5*\spw) -- (-0.45*\sph,0.45*\sph) -- (-0.5*\spw,0.5*\spw) -- cycle;
        \end{scope}
    \end{scope}
}

\newif\ifnomentry

\renewcommand{\nomgroup}[1]{%
  \nomentryfalse
  \ifthenelse{\equal{#1}{V}}{\item[\textbf{Variables}]}{%
    \ifthenelse{\equal{#1}{A}}{\item[\textbf{Abbreviations}]}{}}\nomentrytrue} 
\newcommand{\defvar}[3][show]{\nomenclature[V]{#2}{#3}\ifthenelse{\equal{#1}{show}}{#2}{}}
\newcommand{\defabbrev}[2]{#2 (#1)\nomenclature[A]{#1}{#2}}

\newcommand{\reals}{\mathbb R}

\newcommand{\integers}{\mathbb Z}

\newcommand{\subdiff}[1][none]{%
  \ifthenelse{\equal{#1}{none}}{%
    \partial%
  }{%
    \partial_{#1}%
  }%
}

\newcommand{\diag}[1]{\mathrm{diag}(#1)}

\newcommand{\transp}{{\scriptscriptstyle\mathsf{T}}}

\newcommand{\dd}{\mathrm{d}}

\newcommand{\fun}[2][1]{%
  #2(%
  \foreach \index in {1, ..., #1} {%
    \ifthenelse{\equal{\index}{#1}}{%
      \cdot%
    }{%
      \cdot,%
    }%
  })}

\newcommand{\ev}[3][c]{%
  \ifthenelse{\equal{#1}{c}}{%
    \ifthenelse{\equal{#2}{}}{\forall[0,#3]}{\forall[#2,#3]}%
  }{%
    \ifthenelse{\equal{#1}{o}}{%
      \ifthenelse{\equal{#2}{}}{\forall(0,#3)}{\forall(#2,#3)}%
    }{%
      \ifthenelse{\equal{#1}{oc}}{%
        \ifthenelse{\equal{#2}{}}{\forall(0,#3]}{\forall(#2,#3]}%
      }{%
        \ifthenelse{\equal{#2}{}}{\forall[0,#3)}{\forall[#2,#3)}%
      }%
    }%
  }%
}

\DeclareMathOperator*{\minimize}{\mathrm{minimize}}

\newcounter{l}
\newcounter{j}
\newcounter{k}

\newcommand{\lloptimization}[5][center]{
  \setsepchar{\\}%
  \readlist\mylist{#5}
  \begin{subequations}
    \ifthenelse{\equal{#2}{}}{}{\label{eq:#2}}
    \ifthenelse{\equal{#1}{left}}{
      \begin{flalign}
        \ifthenelse{\equal{#2}{}}{}{\label{eq:#2_a}}
        &\min_{\ifthenelse{\equal{#3}{}}{}{#3}}~ #4\textnormal{ s.t.}\hspace{-1mm} &&\\
        \forloop{k}{0}{\arabic{k} < \listlen\mylist[]}{
          \setcounter{j}{\value{k}+1}
          \setcounter{l}{\value{k}+2}
          \ifthenelse{\equal{#2}{}}{}{\label{eq:#2_\alph{l}}}
          \ifthenelse{\equal{\arabic{j}}{\listlen\mylist[]}}{%
            &\mylist[\arabic{j}]%
          }{%
            &\mylist[\arabic{j}] &&\\%
          }%
        }
      \end{flalign}
    }{
      \begin{align}
        \ifthenelse{\equal{#2}{}}{}{\label{eq:#2_a}}
        \minimize_{\ifthenelse{\equal{#3}{}}{}{#3}}~& #4 \\
        \mathrm{subject~to:~}
        \forloop{k}{0}{\arabic{k} < \listlen\mylist[]}{
        \setcounter{j}{\value{k}+1}
        \setcounter{l}{\value{k}+2}
        \ifthenelse{\equal{#2}{}}{}{\label{eq:#2_\alph{l}}}
        \ifthenelse{\equal{\arabic{j}}{\listlen\mylist[]}}{%
                                                    &\mylist[\arabic{j}]%
                                                      }{%
                                                    &\mylist[\arabic{j}]\\%
        }%
        }
      \end{align}
    }
  \end{subequations}
}

\newcommand{\llpoptimization}[6][left]{
  \begin{problem}\ifthenelse{\equal{#2}{NoNe}}{}{\textbf{#2.}}%
    \ifthenelse{\equal{#3}{}}{}{\label{problem:#3}}%
    \lloptimization[#1]{#3}{#4}{#5}{#6}%
  \end{problem}%
}

\NewEnviron{optimization}[4]{\lloptimization[#1]{#2}{#3}{#4}{\BODY}}
\NewEnviron{foptimization}[5][NoNe]{%
  \begin{mdframed}[default,linewidth=1pt]%
    \llpoptimization[#2]{#1}{#3}{#4}{#5}{\BODY}%
  \end{mdframed}%
}
\NewEnviron{poptimization}[5][NoNe]{%
  \llpoptimization[#2]{#1}{#3}{#4}{#5}{\BODY}%
}

\definecolor{darkolivegreen}{rgb}{0.33, 0.6, 0.18}

\newcommand{\postmeeting}[2][show]{%
  \ifthenelse{\equal{#1}{show}}{%
    {\color{orange} #2}%
  }{}%
}
\newcommand{\idea}[2][show]{%
  \ifthenelse{\equal{#1}{show}}{%
    {\color{orange} #2}%
  }{}%
}
\newcommand{\todo}[2][show]{%
  \ifthenelse{\equal{#1}{show}}{%
    {\color{blue} (TODO: #2)}%
  }{}%
}
\newcommand{\question}[2][show]{%
  \ifthenelse{\equal{#1}{show}}{%
    {\color{darkolivegreen} (Q: #2)}%
  }{}%
}
\newcommand{\fixme}[2][show]{%
  \ifthenelse{\equal{#1}{show}}{%
    {\color{red} (FIXME: #2)}%
  }{}%
}
\newcommand{\mycomment}[2][show]{%
  \ifthenelse{\equal{#1}{show}}{%
    {\color{red} (C: #2)}%
  }{}%
}

\newcommand{\Behcet}{Beh\c{c}et\xspace}
\newcommand{\Acikmese}{A\c{c}{\i}kme\c{s}e\xspace}

\setlist{itemsep=0pt}

\hypersetup{
  colorlinks,
  linkcolor={red},
  citecolor={green},
  urlcolor={blue}
}

\allowdisplaybreaks

\tcbset{highlight math style={enhanced,
colframe=red!60!black,colback=yellow!50!white,arc=4pt,boxrule=1pt,
drop fuzzy shadow}}

\newcolumntype{C}{>{\centering\arraybackslash}X}
\newcolumntype{L}{>{\raggedright\arraybackslash}X}
\newcolumntype{R}{>{\raggedleft\arraybackslash}X}

\newtheorem{theorem}{Theorem}

\theoremstyle{definition}

\newtheorem{definition}{Definition}

\newtheorem{problem}{Problem}

\theoremstyle{empty}

\newmdenv[
outerlinewidth = 1,%
linewidth = 0pt,%
roundcorner = 2pt,%
leftmargin = 20,%
rightmargin = 0,%
backgroundcolor = lightgray!10,%
outerlinecolor = gray!50,%
innertopmargin = \topskip,%
splittopskip = \topskip,%
frametitle = Summary,%
frametitlebelowskip = 0pt,%
]{summary_box}

\newmdenv[
outerlinewidth = 1,%
linewidth = 0pt,%
roundcorner = 2pt,%
leftmargin = 60,%
rightmargin = 60,%
backgroundcolor = yellow!40,%
outerlinecolor = black!50,%
innertopmargin = \topskip,%
splittopskip = \topskip,%
]{highlight_box}

\graphicspath{{./figures/}}

\title{Fast Trajectory Optimization via Successive Convexification for
  Spacecraft Rendezvous with Integer Constraints}

\author{
  Danylo Malyuta\footnote{Doctoral Student, W.E. Boeing Department of
    Aeronautics \& Astronautics \texttt{danylo@uw.edu}}, %
  Taylor P. Reynolds\footnote{Doctoral Student, W.E. Boeing Department of
    Aeronautics \& Astronautics \texttt{tpr6@uw.edu}}, %
  Michael Szmuk\footnote{Doctoral Student, W.E. Boeing Department of Aeronautics
    \& Astronautics \texttt{mszmuk@uw.edu}}, %
  \Behcet \Acikmese\footnote{Professor, W.E. Boeing Department of Aeronautics \&
    Astronautics, AIAA Associate Fellow, \texttt{behcet@uw.edu}} and %
  Mehran Mesbahi\footnote{Professor, W.E. Boeing Department of Aeronautics \&
    Astronautics, AIAA Associate Fellow, \texttt{mesbahi@uw.edu}}
}
\affil{University of Washington, Seattle, WA, 98195}

\makenomenclature

\begin{document}

\maketitle

\begin{abstract}
  In this paper we present a fast method based on successive convexification for
  generating fuel-optimized spacecraft rendezvous trajectories in the presence
  of mixed-integer constraints. A recently developed paradigm of state-triggered
  constraints allows to efficiently embed a subset of discrete decision
  constraints into the continuous optimization framework of successive
  convexification. 
  As a result, we are able to solve difficult trajectory optimization problems
  at interactive speeds, as opposed to a mixed-integer programming approach that
  would require significantly more solution time and computing power. Our method
  is applied to the real problem of transposition and docking of the Apollo
  command and service module with the lunar module. We demonstrate that, within
  seconds, we are able to obtain trajectories that are up to 90~percent more
  fuel efficient (saving up to 45~kg of fuel) than non-optimization based
  Apollo-era design targets. Our trajectories take explicit account of minimum
  thrust pulse width and plume impingement constraints. Both of these
  constraints are naturally mixed-integer, but we handle them as state-triggered
  constraints. In its current state, our algorithm will serve as a useful
  off-line design tool for rapid trajectory trade studies.
\end{abstract}

\begin{multicols}{2}
\printnomenclature
\end{multicols}



\section{Introduction}
\label{sec:introduction}


\lettrine{S}{pace} programs have historically been deemed mature once they establish the ability to perform rendezvous and docking operations~\cite{Woffinden2007}. Some of the earliest programs of the United States and Soviet Union (e.g., Gemini, Soyuz) had as their explicit goal to demonstrate the capability of performing rendezvous, proximity operations, and docking maneuvers. The ultimate objective to land humans on the moon drove the need for these capabilities. Beyond the lunar missions of the 1960s, rendezvous and docking continued to be a core technology required to construct and service space stations that were built in low Earth orbit~\cite{Goodman2006}. The Shuttle program was comprised of dozens of missions for which rendezvous (and more generally, proximity operations) was an explicit mission objective. The core technology used to achieve rendezvous and docking has remained largely unchanged in the decades since the earliest maneuvers were successful. While this heritage technology is far from obsolete, it has been stated that it may be unable to meet the requirements of future missions~\cite{Woffinden2007}. A
driving force that will require new methods is the need for a system that can
perform \textit{fully autonomous} rendezvous in several domains (e.g., low Earth
orbit, low lunar orbit, etc.)~\cite{DSouza2007}.

The objective of this paper is to present a framework for designing autonomous
docking trajectories that accurately reflect the capabilities and constraints
that have been historically prevalent for proximity operation missions. We view
the problem as a trajectory generation problem, and compute what would be implemented as guidance solutions. Explicit consideration of the navigation and closed-loop control aspects of the problem are beyond the scope of this paper. We show how to model challenging constraints within a continuous optimization framework, and the result is an algorithm that can be used as a design tool while performing trade studies for autonomous rendezvous maneuvers. 

The open-loop generation of spacecraft docking trajectories using
optimization-based methods has been studied only recently, spawned naturally by
the shift towards autonomy. Open-loop trajectory generation is a slightly
different (but intimately related) technique than receding horizon or model
predictive control, which has received more attention in the rendezvous and
docking
literature~\cite{Weiss2014,Richards2002,Hartley2013}. In~\cite{Miele2007,Miele2007a}
the authors discussed both time- and fuel-optimal solutions with a focus on
problem formulations that were conducive to on-board implementation. Their study
offers an insightful view on the structure of optimality at the cost of a
simplified problem statement and omission of state constraints. In
\cite{Pascucci2017}, lossless convexification is used to generate fuel-optimal
docking trajectories which account for non-convex thrust and plume impingement
constraints, albeit the thrust is not allowed to turn off. In \cite{Breger2008},
an optimization framework is used to impose safety-based constraints in the case
of anomalous behaviour (including thruster failure) by introducing a suboptimal
convex program to design safe trajectories which approximate a non-convex
mixed-integer problem using a new set of ``safe'' inputs. Along the same lines
of mixed-integer programming,~\cite{Richards2002} solved a fuel-optimal problem
subject to thrust plume and collision avoidance constraints. They introduced
several heuristic techniques in order to fit the problem within the scope of
mixed-integer linear programming, but still observed rather large solve times
(over 40 minutes in some cases). More recently,~\cite{Sun2019} studied a
multi-phase docking problem with several state constraints. The authors used
binary variables to impose different constraints during each phase, and proposed
an iterative solution method with closed-form update rules. Beyond the use of
mixed-integer methods,~\cite{Phillips2003} proposed a randomized optimization
method similar to the $A^*$ method, while~\cite{Hartley2013} proposed a convex
$\ell_1$-regularized model predictive control solution.

Notably, each of the aforementioned references do not (i) consider the
spacecraft attitude during trajectory generation and (ii) explicitly account for
what is referred to as the \defabbrev{MIB}{minimum impulse bit} of the reaction
control thrusters that are used to realize the trajectories. The latter
constraint refers to the fact that these impulsive thrusters cannot realize an
arbitrarily low thrust; there is some minimum pulse-width that is inherent to
the hardware. Ref.~\cite{Hartley2013} acknowledges this issue, but uses an
$\ell_1$ penalty term to discourage violation of the constraint (i.e., a soft
constraint), instead of explicitly enforcing it.

\subsection{Contributions}
\label{subsec:contributions}

The contribution of this paper is a numerical solution of the 6-\defabbrev{DoF}{degree-of-freedom} autonomous rendezvous and docking problem
with consideration for the MIB of the \defabbrev{RCS}{reaction control system},
plume impingement constraints, and several state constraints that are common to
docking trajectories (e.g., the approach cone). Moreover, we demonstrate
the use of continuous optimization techniques (no integer variables) that
can be solved fast enough to facilitate rapid trade studies and trajectory
design. By using the newly introduced state-triggered
constraints~\cite{Szmuk2018,Reynolds2019b}, we show how to model the MIB
constraint within a continuous optimization problem, effectively solving the
problem identified (but not solved) in~\cite{Hartley2013}. Moreover, plume
impingement constraints that restrict the directions in which thrusters may be
fired are enforced. These important constraints were a primary driver for the
design of the Shuttle docking maneuvers~\cite{Woffinden2007,Goodman2006}, but it
has been noted that they are not explicitly necessary unless the two vehicles
are in close proximity~\cite{Richards2002}. As such, we again use
state-triggered constraints to enforce plume impingement constraints
\textit{only} when the two vehicles are close enough to one another, and stress that this is achieved without the use of mixed-integer programming or multi-phase optimization.

The sensor suite that is used to perform vehicle navigation during autonomous
rendezvous may require specific pointing requirements based on the distance
between the docking vehicle and target. For example, when optical sensing is
used for the terminal docking phase~\cite{Woffinden2007,Fehse2003}, the camera
must be pointed towards the docking target. In~\cite{Fehse2003} it is stated
explicitly that ``the nominal attitude of the chaser vehicle is determined by
[...] the operational range of the sensors for attitude and trajectory control
[and] by the range of the antennas for communication with ground and with the
target station''. Using state-triggered constraints, it is readily possible to
constrain the set of chaser feasible attitudes based on navigation sensor range
without recursion to binary variables or multi-phase optimization. The plume
impingement constraint in this paper serves to demonstrate this approach.

\subsection{Outline}
\label{subsec:outline}

The outline of this paper is as follows. First,
Section~\ref{sec:problem_formulation} details the rendezvous problem
formulation, including the 6-DoF kinematics and dynamics, the impulsive thrust
model, and the state constraints that are imposed. A statement of the
fixed-final time (non-convex) optimal control problem is provided in
Section~\ref{subsec:minlp}. Section~\ref{sec:solution_method} outlines our
solution method based on successive convexification for solving said optimal
control problem. The analytic expressions that are used to obtain solutions that
respect the non-linear dynamics of the problem are provided in
Section~\ref{subsec:propagation}. The use of state-triggered constraints for
dealing with discrete decision in a continuous optimization framework is also
explained in Section~\ref{subsec:stc}. A statement of the convexified
sub-problem is given in Section~\ref{subsec:cvx}. Lastly,
Section~\ref{sec:example} corroborates the effectiveness of our algorithm by
using an example of the Apollo command and service module transposition and
docking maneuver with the lunar module. Within a few seconds, as a consequence
of the optimization-based approach, our method yields solutions that are up to
$90~\%$ more fuel efficient than the Apollo design target.

\subsection{Notation}
\label{subsec:notation}

We use the following notation and conventions. Let $\reals$, $\integer$ and
$\quaternion$ denote the sets of reals, integers and unit quaternions
respectively. We use $\reals_{\ge 0}$ and $\reals_{>0}$ for the non-negative and
the positive reals, and similarly for the integers. Scalars and vectors are
lowercase, e.g. $x\in\reals$ or $y\in\reals^n$, matrices are uppercase,
e.g. $M\in\reals^{n\times m}$, and coordinate frames are calligraphic subscripts
of $\mathcal F$, e.g. $\fr{X}$ is the ``$\mathcal X$ frame''. We generally use
the following convention for subscripts and superscripts: ${}_{j}x_k^i$ denotes
the $i$-th element of variable $x$ at the $k$-th discrete time step and at the
$j$-th iteration of the successive convexification algorithm. When an index is
\textit{unspecified}, e.g. $k$ is not explicitly assigned a set of values, then
the notation denotes a trajectory such that, based on context, ${}_{j}x_k^i$ may
refer to the discrete state trajectory
$\{{}_{j}x_0^i,{}_{j}x_1^i,{}_{j}x_2^i,\dots\}$. The operators $\diag{\cdot}$
and $\blkdiag{\cdot,\cdot,\dots}$ build diagonal and block-diagonal matrices in
the same way as they do in high-level programming languages such as MATLAB and
Python's Numpy. The vector $e_i\in\reals^n$ is the $i$-th Eucledian basis unit
vector whose $i$-th element is one and where all other elements are zero. The
vector $\bm{1}_n\in\reals^n$ is the $n$-dimensional vector of ones,
$I_n\in\reals^{n\times n}$ is the identity matrix of size $n$, and
$0_{n\times m}\in\reals^{n\times m}$ is the zero matrix or vector. In all
figures, we use the red/green/blue colors to denote the $+x$/$+y$/$+z$ axes or
vector components.

\section{Rendezvous Problem Formulation}
\label{sec:problem_formulation}

In this section we formulate the problem of guiding a dynamic \textit{chaser}
spacecraft to dock with a passive \textit{target} spacecraft whose trajectory is
predetermined. We define the chaser's dynamics in Section~\ref{subsec:dynamics},
its actuator model in Section~\ref{subsec:input_model}, and the rendezvous
constraints in Sections~\ref{subsec:impingement} and
\ref{subsec:approach_cone}. Section~\ref{subsec:minlp} gives a complete
formulation of the fixed-final time non-convex optimal control problem which, if
solved, generates the fuel-optimal rendezvous trajectory.

\subsection{Chaser Spacecraft Dynamics}
\label{subsec:dynamics}

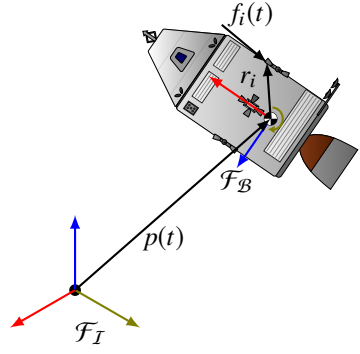
\begin{figure}[tb]
  \centering
    \begin{tikzpicture}
\tikzmath{
	\cx=-2; \cy=-2; \wdth=8; \hght=8;
	\scx=1; \scy=0; \tlt=55.5;
    \comscl=0.6; \comx=\scx-0.4; \comy=\scy+0.28;
    \len=1; \rot=150;
    \pX=  0; \Xx=cos(\pX)*\len; \Xy=sin(\pX))*\len;
    \pY=120; \Yx=cos(\pY)*\len; \Yy=sin(\pY))*\len;
    \pZ=240; \Zx=cos(\pZ)*\len; \Zy=sin(\pZ))*\len;
}


\isoaxesNL{\cx}{\cy}{210}{1}{0}{$\fr{I}$}

\apolloCSM{\scx}{\scy}{\tlt}{0}{0.4}

    \begin{scope}[shift={(\comx,\comy)}]
        \draw[Xcol,thick,->,rotate=\tlt] (0,0) -- +(0,1);
        \draw[Ycol,thick,{Latex[scale=0.4,angle'=60,flex=1]}-] (-0.04,-0.4em) arc (-110:85:0.15);
        \draw[Zcol,thick,->,rotate=\tlt+90] (0,0) -- +(0,0.8);
        \begin{scope}[scale=\comscl]
            \filldraw[black,radius=0.4em] (0,0) circle;%
            \fill[fill=white,radius=0.4em] (0,0) -- ++(0.4em,0) arc [start angle=0,end angle=90] -- ++(0,-0.8em) arc [start angle=270, end angle=180];%
        \end{scope}
    \end{scope}
\node[name=txt] at (\scx-0.85,\scy-0.55) {$\fr{B}$};

\coordinate (RCS) at (\comx-0.05,\comy+0.75);
\draw[thick,->,black] (\comx,\comy) -- (RCS) node[anchor=east,pos=0.7] {$r_i$};
\draw[thick,<-,black] (RCS) -- ++(-0.6,0.5) node[anchor=south,pos=0.7,xshift=0.2cm] {$f_i(t)$};
\draw[thick,->] (-2,-2) -- (\comx,\comy) node[anchor=west,pos=0.3] {$p(t)$};


\end{tikzpicture}
  \caption{Illustration of the inertia and body frames, and the $i$-th RCS
    thruster.  The thruster is represented by its position $r_i$ relative to the
    CoM and its thrust vector $f_i(t)$, both quantities expressed in the body
    frame.}
  \label{fig:general_thruster_geometry}
\end{figure}

We now develop the first-order differential equations that govern the chaser
spacecraft's position and attitude. Let \defvar{$\fr{B}$}{chaser body frame} be
a body-fixed frame centered at the chaser's \defabbrev{CoM}{center of
  mass}. Assume that the chaser has a constant mass $m\in\reals_{>0}$ and an
inertia tensor $J\in\reals^{3\times 3}$. We model the position dynamics as a
double-integrator for simplicity. We note that more refined translational
dynamics (e.g. Hill's equations) can be incorporated into the design easily
without changing the solution methodology. The attitude dynamics are modeled by
Euler's equations, yielding the overall dynamics: \defvar[hide]{$m$}{chaser
  mass} \defvar[hide]{$J$}{chaser inertia tensor}
\begin{subequations}
  \label{eq:vanilla_dynamics}
  \begin{align}
    \label{eq:dpdt}
    \dot p(t) &= v(t), \\
    \label{eq:dvdt}
    \dot v(t) &= \frac{1}{m}\sum_{i=1}^Mq(t)\otimes f_i(t)\otimes q(t)^*, \\
    \label{eq:dqdt}
    \dot q(t) &= \frac{1}{2}q(t)\otimes\omega(t), \\
    \label{eq:dwdt}
    \dot\omega(t) &= J^{-1}\biggl[\sum_{i=1}^M r_i\times f_i(t)-
                    \omega(t)\times (J\omega(t))\biggr],
  \end{align}
\end{subequations}
where $q(t)\in\quaternion$ is the quaternion representation of a frame change
from $\fr{B}$ to the inertial frame \defvar{$\fr{I}$}{inertial frame}, and
$\otimes$ is quaternion multiplication. We use the Hamilton quaternion
convention \cite{Sola2017}. Position and attitude control occurs using
\defvar{$M$}{number of RCS thrusters} RCS thrusters whose operation is described
in Section~\ref{subsec:input_model}. As illustrated in
Figure~\ref{fig:general_thruster_geometry}, the $i$-th thruster is located at
position $r_i\in\reals^3$ in the body frame, and produces a thrust
\defvar{$f_i(t)$}{thrust vector of $i$-th RCS thruster in $\fr{B}$} and a torque
$r_i \times f_i(t)$ of variable duration, each expressed in
$\fr{B}$. \defvar[hide]{$r_i$}{position of $i$-th RCS thruster in $\fr{B}$}

\subsection{Impulsive Thrust Model}
\label{subsec:input_model}

\begin{figure}[tb]
  \centering
  \includegraphics[width=0.6\textwidth]{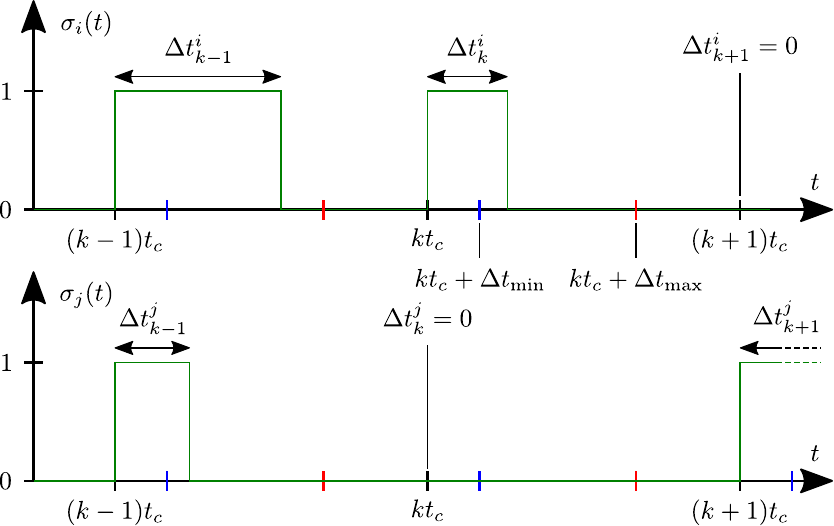}
  \caption{Illustration of the pulsed thrust model, with the $i$-th and $j$-th
    thrusters shown. Each thruster can be off or produce a thrust vector
    $\hat f_i$ for at least $\Delta t_{\min}$ and at most $\Delta t_{\max}$
    seconds. Pulses are spaced $t_c$ seconds apart.}
  \label{fig:pulse_model}
\end{figure}

Each RCS thruster is assumed to be capable of producing short constant-thrust
pulses along a body-fixed direction vector. The $i$-th thrust can thus be
written as
\begin{equation}
  \label{eq:thrust_model}
  f_i(t) = \sigma_i(t)\hat f_i, \quad \sigma_i(t)\definedas
  \begin{cases}
    1 & \text{if }t\in [k t_c,k t_c+\Delta t_k^i]~\text{for some}~k\in\integer, \\
    0 & \text{otherwise,}
  \end{cases}
\end{equation}
where \defvar{$t_c$}{control interval duration} is a \textit{control interval}
which defines the duration between consecutive pulses and
\defvar{$\Delta t_k^i$}{pulse width of thruster $i$ for control interval $k$} is
the pulse width of the $i$-th thruster for the $k$-th control
interval. Figure~\ref{fig:pulse_model} illustrates the concept. Due to delays in
on-board electronics and residual propellant flow downstream of the injector
valves \cite[pp.~2.5-16~to~2.5-18]{aoh_vol1}, the pulse width is lower bounded
such that $\Delta t_k^i\ge$\defvar{$\Delta t_{\min}$}{minimum pulse width} if
$\Delta t_k^i\ne 0$. Other propulsion and RCS parameters, such as engine service
life and damage to engine materials, may impose an upper bound
$\Delta t_k^i\le$\defvar{$\Delta t_{\max}$}{maximum pulse width}. As a result,
the following constraint must be satisfied:
\begin{equation}
  \label{eq:pulse_width_constraint}
  \begin{rcases}
  \Delta t_k^i<\Delta t_{\min}\Rightarrow \Delta t_k^i = 0\ \\
  0 \leq \Delta t_k^i \leq \Delta t_{\max} 
  \end{rcases}\quad\text{for all}~i=1,\dots,M\text{ and }k\in\integer_{\ge 0}.
\end{equation}

\subsection{Plume Impingement Constraint}
\label{subsec:impingement}

A plume impingement constraint prevents the chaser's RCS thrusters from firing
at the target spacecraft. As mentioned previously, this constraint is only required when the chaser is close enough to the target for impingement to be of concern. Let $q_f\in\quaternion$\defvar[hide]{$(p_f,v_f,q_f,\omega_f)$}{final
  chaser state} be the desired chaser final attitude and assume that we want to
activate the impingement constraint when the chaser is within an
$r_a\in\reals$\defvar[hide]{$r_a$}{approach radius} \textit{approach radius}
away from its final docked position $p_f\in\reals^3$. We use the error between the chaser's attitude and the desired final docked attitude as a proxy for the impingement constraint. The error quaternion between the current and final attitudes is given by
$\Delta q(t)=q(t)^*\otimes q_f\in\quaternion$. Denoting by $[q_f]_\otimes$ the
right quaternion product matrix, the impingement constraint is written as:
\begin{equation}
  \label{eq:impingement_attitude}
  \|p(t)-p_f\|_2< r_a\Rightarrow
  e_1^\transp [q_f]_\otimes q(t)^*\ge\cos(\Delta\theta_{\max}/2),
\end{equation}
where \defvar{$\Delta\theta_{\max}$}{impingement maximum attitude error}
is the maximum error angle that the chaser can have with respect to its terminal
attitude when it is within the approach radius. Suppose also that there is a
subset $\mathcal M\in\{1,\dots,M\}$\defvar[hide]{$\mathcal M$}{forward-facing
  thruster indices} of forward-facing RCS thrusters that are to be kept silent
once the chaser is within the approach radius of the target. We impose this
constraint as:
\begin{equation}
  \label{eq:impingement_forward_thrusters}
  \|p(t)-p_f\|_2< r_a\Rightarrow\sigma_i(t)=0, \quad \text{for all}~i\in\mathcal M.
\end{equation}

\subsection{Approach Cone Constraint}
\label{subsec:approach_cone}

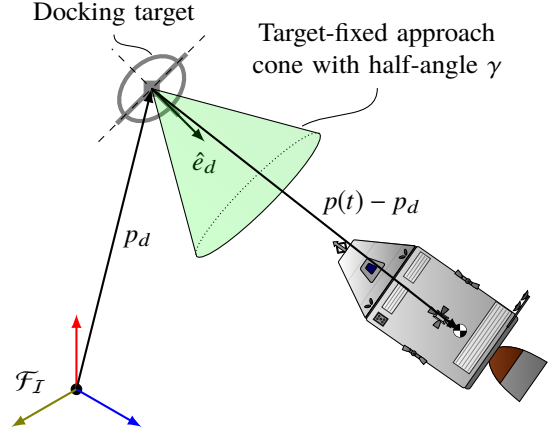
\begin{figure}[tb]
    \centering
    \begin{tikzpicture}
\tikzmath{
	\cx=-3; \cy=-2; \wdth=8; \hght=8;
	\scx=2.5; \scy=-1.5; \tlt=55;
	\tx=-2; \ty=2;
	\magn=sqrt((\tx-\scx)*(\tx-\scx)+(\ty-\scy)*(\ty-\scy));
	\scux=(\scx-\tx)/\magn; \scuy=(\scy-\ty)/\magn;
    \Acone=45; \Lcone=2; \Dcone=0.2; \apang=60;
    \edx=cos(\Acone); \edy=-sin(\Acone);
    \comscl=0.6; \comx=\scx-0.4; \comy=\scy+0.28;
}


\isoaxesNL{\cx}{\cy}{90}{1}{0}{$\fr{I}$}

\apolloCSM{\scx}{\scy}{\tlt}{0}{0.4}
\begin{scope}[shift={(\comx,\comy)},scale=\comscl]
    \filldraw[black,radius=0.4em] (0,0) circle;%
    \fill[fill=white,radius=0.4em] (0,0) -- ++(0.4em,0) arc [start angle=0,end angle=90] -- ++(0,-0.8em) arc [start angle=270, end angle=180];%
\end{scope}

\begin{scope}[shift={(\tx,\ty)},rotate=45]
\draw[color=white!50!black,ultra thick] (0,0) circle (0.5 and 0.35); 
\filldraw[color=white!50!black] (-0.65,-0.025) rectangle (0.65,0.025);
\filldraw[color=white!50!black] (0.025,0) rectangle (-0.025,-0.5);
\filldraw[color=white!50!black,rotate=45] (-0.1,-0.1) rectangle (0.1,0.1); 
\draw[dashed] (-1,0) -- (1,0);
\draw[dashed] (0,-0.75) -- (0,0.75);
\end{scope}

\draw[thick,->] (-3,-2) -- (\tx,\ty) node [anchor=west,pos=0.5] {$p_d$};
\draw[thick,->] (\tx,\ty) -- ++(\edx,\edy) node [anchor=north] {$\hat{e}_d$};
\cone{\tx}{\ty}{\Acone}{\Lcone}{\apang}{\Dcone}{color=green!100,opacity=0.10,shading=axis,shading angle=90,left color=green!100,right color=green!70!black}
\coneback{\tx}{\ty}{\Acone}{\Lcone}{\apang}{\Dcone}{color=black!100,densely dotted}
\draw[thick,->] (\tx,\ty) -- (\comx,\comy) node [anchor=west,pos=0.5,xshift=0.1cm,yshift=0.1cm] {$p(t)-p_d$};

\node[align=center,text width=4cm,name=txt] at (\tx+3,\ty+0.5) {Target-fixed approach cone with half-angle $\gamma$};
\draw[] (txt.south) to[out=-90,in=30] (\tx+2,\ty-0.5);
\node[align=center,text width=4cm,name=dtrg] at (\tx-0.5,\ty+1) {Docking target};
\draw[] (dtrg.south) to[out=-90,in=125] (\tx,\ty+0.45);

\end{tikzpicture}
    \caption{An approach cone emanating from the docking port constrains the positions of the chaser.}
    \label{fig:approach_cone}
\end{figure}

The approach cone constraint ensures that the chaser is sufficiently in front of
the target to ensure that a line of sight exists between the chaser's sensors
and the docking target. As shown in Figure~\ref{fig:approach_cone}, let
$p_d\in\reals^3$\defvar[hide]{$(p_l,v_l,q_l,\omega_l)$}{lunar module state} be
the docking port location, $\hat e_d\in\reals^3$
\defvar[hide]{$\hat e_d$}{docking axis in $\fr{I}$} the docking axis in
$\fr{I}$, and $\gamma\in (0,\pi/2)$\defvar[hide]{$\gamma$}{approach cone
  half-angle} be the approach cone half-angle. The approach cone constraint is
written as:
\begin{equation}
  \label{eq:approach_cone}
  \|p(t)-p_d\|_2\cos(\gamma)\le (p(t)-p_d)^\transp\hat e_d.
\end{equation}

\subsection{Basic Rendezvous Guidance Problem}
\label{subsec:minlp}

We plan a rendezvous trajectory of a fixed duration \defvar{$t_f$}{rendezvous
  duration}\defvar[hide]{$(p_0,v_0,q_0,\omega_0)$}{initial chaser state} and fixed boundary conditions:
\begin{gather}
  \label{eq:initial_condition}
  p(0)=p_0,\quad v(0)=v_0,\quad q(0)=q_0,\quad \omega(0)=\omega_0, \\
  p(t_f)=p_f,\quad v(t_f)=v_f,\quad q(t_f)=q_f,\quad \omega(t_f)=\omega_f.
\end{gather}

Although the duration $t_f$ is fixed, an optimal value of $t_f$ can be searched
for via a line search method. A free-final time formulation can also be used,
which has been studied in previous work on successive convexification
\cite{Szmuk2018,Reynolds2019b}. We use a cost which minimizes the cumulative
thruster firing time and is a proxy for total fuel usage
\cite[Section~4.3.4.1.2]{csm_aoh}: \defvar[hide]{$J_f$}{minimum pulse width
  cost}
\begin{equation}
  \label{eq:original_cost}
  J_f = \sum_{i=1}^M\int_0^{t_f}\sigma_i(t)\dd t.
\end{equation}

We are now able to write a continuous-time, fixed-final time non-convex optimal
control problem which we call the \defabbrev{BRGP}{basic rendezvous trajectory
  problem}:

\begin{poptimization}[Basic Rendezvous Guidance
  Problem]{center}{brtop}{\sigma_i(t)}{J_f}
  \dot p(t) = v(t), \\
  \dot v(t) = \frac{1}{m}\sum_{i=1}^Mq(t)\otimes f_i(t)\otimes q(t)^*, \\
  \dot q(t) = \frac{1}{2}q(t)\otimes\omega(t), \\
  \dot\omega(t) = J^{-1}\biggl[\sum_{i=1}^M r_i\times f_i(t)-
  \omega(t)\times (J\omega(t))\biggr], \\
  0\le\Delta t_k^i\le\Delta t_{\max}~\text{for
    all}~i=1,\dots,M\text{ and }k\in\integer_{\ge 0}, \\
  \Delta t_k^i<\Delta t_{\min}\Rightarrow \Delta t_k^i=0~\text{for
    all}~i=1,\dots,M\text{ and }k\in\integer_{\ge 0}, \\
  \|p(t)-p_f\|_2< r_a\Rightarrow
  e_1^\transp [q_f]_\otimes q(t)^*\ge\cos(\Delta\theta_{\max}/2), \\
  \|p(t)-p_f\|_2< r_a\Rightarrow\sigma_i(t)=0~\text{for all}~i\in\mathcal M, \\
  \|p(t)-p_d\|_2\cos(\gamma)\le (p(t)-p_d)^\transp\hat e_d, \\
  p(0)=p_0,~v(0)=v_0,~q(0)=q_0,~\omega(0)=\omega_0, \\
  p(t_f)=p_f,~v(t_f)=v_f,~q(t_f)=q_f,~\omega(t_f)=\omega_f.
\end{poptimization}

\section{Solution via Successive Convexification}
\label{sec:solution_method}

In this section we propose an algorithm for finding a feasible solution of
Problem~\ref{problem:brtop} which also attempts to locally minimize $J_f$ in
\eqref{eq:original_cost}. 
Sections~\ref{subsec:propagation}-\ref{subsec:cvx} formulate a local
approximation of Problem~\ref{problem:brtop} as a convex optimization
sub-problem. This sub-problem then gets used in an iterative solution scheme
called successive convexification, which is detailed in
Section~\ref{subsec:ftr}.

\subsection{Dynamics Linearization and Discretization}
\label{subsec:propagation}

Two obstacles that prevent a fast numerical solution of
Problem~\ref{problem:brtop} are rooted in the nature of the dynamics
\eqref{eq:vanilla_dynamics}:
\begin{enumerate}
\item The dynamics are non-convex due to a rotation operation in \eqref{eq:dvdt}
  and \eqref{eq:dqdt}, and the cross product term in \eqref{eq:dwdt};
\item The dynamics evolve in continuous-time, thus the solution space is infinite
  dimensional. Numerical optimization algorithms can only handle problems where
  the number of decision variables is finite.
\end{enumerate}

The first obstacle is resolved by linearizing, and the second obstacle is
resolve by discretizing the dynamics. We first discuss linearization since it is
a necessary precursor to discretization. To begin, define the following
functions based on the dynamics \eqref{eq:vanilla_dynamics}:
\begin{align}
  \label{eq:f_omega}
  f_\omega(\omega(t)) &\definedas -J^{-1}[\omega(t)\times (J\omega(t))], \\
  \label{eq:f_q}
  f_q(q(t),\omega(t)) &\definedas \frac{1}{2}q(t)\otimes\omega(t), \\
  \label{eq:f_v}
  f_v(q(t),f_i(t)) &\definedas \frac{1}{m}\sum_{i=1}^Mq(t)\otimes f_i(t)\otimes q(t)^*.
\end{align}

Let
$(\bar p(\cdot),\bar v(\cdot),\bar q(\cdot),\bar\omega(\cdot)):[0,t_f]\to\reals^{13}$
be a given reference state trajectory and $\bar f_i(\cdot):[0,t_f]\to\reals^{3}$
be a given reference input trajectory. Section~\ref{subsec:ftr} will explain how
these reference trajectories are obtained. The Jacobians of
\eqref{eq:f_omega}-\eqref{eq:f_v} are:
\begin{equation}
  A_{\omega\omega}(t) = \left.\frac{\partial f_\omega}{\partial\omega}\right|_{\bar\omega(t)},\quad
  A_{qq}(t) = \left.\frac{\partial f_q}{\partial q}\right|_{\bar q(t),\bar\omega(t)},\quad
  A_{q\omega}(t) = \left.\frac{\partial f_q}{\partial\omega}\right|_{\bar q(t),\bar\omega(t)},\quad
  A_{vq}(t) = \left.\frac{\partial f_v}{\partial q}\right|_{\bar q(t),\bar f_i(t)}.
\end{equation}

As a result, the linearized dynamics in differential form can be written as:
\begin{subequations}
  \label{eq:linearized_differential_form_dynamics}
  \begin{alignat}{3}
    \label{eq:dpdt_linear_differential}
    \dot p(t) &= v(t), \\
    \label{eq:dvdt_linear_differential}
    \dot v(t) &= A_{vq}(t)q(t)+\frac{1}{m}\sum_{i=1}^M\bar q(t)\otimes f_i(t)\otimes
    \bar q(t)^*+r_v(t),~&&r_v(t)\definedas -A_{vq}(t)\bar q(t), \\
    \dot q(t) &= A_{qq}(t)q(t)+A_{q\omega}(t)\omega(t)+r_q(t),~
    &&r_q(t)\definedas
    f_q(\bar q(t),\bar\omega(t))-
    A_{qq}(t)\bar q(t)-A_{q\omega}(t)\bar\omega(t), \\
    \dot\omega(t) &= A_{\omega\omega}(t)\omega(t)+J^{-1}\sum_{i=1}^M r_i\times f_i(t)+r_\omega(t),~&&r_\omega(t)\definedas
    f_{\omega}(\bar\omega(t))-A_{\omega\omega}(t)\bar\omega(t).
  \end{alignat}
\end{subequations}

Let us consider the $k$-th control interval as $[kt_c,(k+1)t_c]$. Then, the
linearized dynamics \eqref{eq:linearized_differential_form_dynamics} can be
written in integral form by leveraging the state transition matrix
\cite{Malyuta2019,Antsaklis2007}:
\begin{subequations}
  \label{eq:linearized_integral_form_dynamics}
  \begin{align}
    \label{eq:dpdt_linear_integral}
    p(t) &= p(kt_c)+\int_{kt_c}^t v(\tau)\dd\tau, \\
    \label{eq:dvdt_linear_integral}
    v(t) &= v(kt_c)+\int_{kt_c}^t A_{vq}(\tau)q(\tau)\dd\tau+
           \frac{1}{m}\sum_{i=1}^M\int_{kt_c}^t\bar q(\tau)\otimes f_i(\tau)\otimes\bar q(\tau)^*\dd\tau+\int_{kt_c}^t r_v(\tau)\dd\tau, \\
    q(t) &= \Phi_q(t,kt_c)q(kt_c)+\int_{kt_c}^t\Phi_q(t,\tau)A_{q\omega}(\tau)\omega(\tau)\dd\tau+
           \int_{kt_c}^t\Phi_q(t,\tau)r_q(\tau)\dd\tau, \\
    \omega(t) &= \Phi_\omega(t,kt_c)\omega(kt_c)+
                J^{-1}\sum_{i=1}^M\int_{kt_c}^t r_i\times f_i(\tau)d\tau+
                \int_{kt_c}^t\Phi_\omega(t,\tau)r_\omega(\tau)\dd\tau,
  \end{align}
\end{subequations}
where the state transition matrices satisfy the following dynamics:
\begin{subequations}
  \begin{align}
    \dot\Phi_q(t,kt_c) &= A_{qq}(t)\Phi_q(t,kt_c),~\Phi_q(kt_c,kt_c)=I\in\reals^{4\times 4}, \\
    \dot\Phi_{\omega}(t,kt_c) &= A_{\omega\omega}(t)\Phi_{\omega}(t,kt_c),~\Phi_{\omega}(kt_c,kt_c)=I\in\reals^{3\times 3}.
  \end{align}
\end{subequations}

Discretization of the dynamics involves computing a \textit{state update}
equation which gives $(p((k+1)t_c),v((k+1)t_c),q((k+1)t_c),\omega((k+1)t_c))$ as
a function of $(p(kt_c),v(kt_c),q(kt_c),\omega(kt_c))$ and the input choice
$f_i(\cdot):[kt_c,(k+1)t_c]\to\reals^3$. This may be done by evaluating
\eqref{eq:linearized_integral_form_dynamics} at $(k+1)t_c$, i.e. at the
end of the control interval. By explicitly using the thruster model
\eqref{eq:thrust_model}, the quantities in \eqref{eq:linearized_integral_form_dynamics}
evaluated at time $(k+1)t_c$ become:
\begin{subequations}
  \label{eq:linearized_integral_form_dynamics_at_tc}
  \begin{align}
    \label{eq:p_tc}
    p((k+1)t_c) &= p(kt_c)+\int_{kt_c}^{(k+1)t_c} v(\tau)\dd\tau, \\
    \label{eq:v_tc}
    v((k+1)t_c) &= v(kt_c)+\int_{kt_c}^{(k+1)t_c} A_{vq}(\tau)q(\tau)\dd\tau+
             \frac{1}{m}\sum_{i=1}^M\int_{kt_c}^{kt_c+\Delta t_k^i}\bar q(\tau)\otimes
             \bar f_i\otimes\bar q(\tau)^*\dd\tau+\int_{kt_c}^{(k+1)t_c} r_v(\tau)\dd\tau, \\
    \nonumber
    q((k+1)t_c) &= \Phi_q((k+1)t_c,kt_c)q(kt_c)+\int_{kt_c}^{(k+1)t_c}\Phi_q((k+1)t_c,\tau)A_{q\omega}(\tau)\omega(\tau)\dd\tau+ \\
    \label{eq:q_tc}
    &\phantom{= \Phi_q((k+1)t_c,kt_c)q(kt_c)+}
             \int_{kt_c}^{(k+1)t_c}\Phi_q((k+1)t_c,\tau)r_q(\tau)\dd\tau, \\
    \nonumber
    \omega((k+1)t_c) &= \Phi_\omega((k+1)t_c,kt_c)\omega(kt_c)+
                       J^{-1}\sum_{i=1}^M\int_{kt_c}^{kt_c+\Delta t_k^i} \Phi_\omega((k+1)t_c,\tau)(r_i\times\hat f_i)d\tau+ \\
    \label{eq:w_tc}
    &\phantom{= \Phi_\omega((k+1)t_c,kt_c)\omega(kt_c)+}
                \int_{kt_c}^{(k+1)t_c}\Phi_\omega((k+1)t_c,\tau)r_\omega(\tau)\dd\tau,
  \end{align}
\end{subequations}
where we observe from \eqref{eq:v_tc} and \eqref{eq:w_tc} that the effect of the
thrust pulse width (our chosen control variable) is to alter the upper limit of
the integration. Since this effect is non-linear, a further linearization is
necessary with respect to $\Delta t_k^i$. This is done via the Leibniz integral
rule \cite[Theorem~3]{Protter1985}. The result is the desired set of update
equations:
\begin{subequations}
  \label{eq:dltv}
  \begin{align}
    \label{eq:p+}
    p_{k+1} &= p_k+A_{pv,k}v_k+A_{pq,k}q_k+A_{p\omega,k}\omega_k+B_{p,k}u_k+r_{p,k}, \\
    \label{eq:v+}
    v_{k+1} &= v_k+A_{vq,k}q_k+A_{v\omega,k}\omega_k+B_{v,k}u_k+r_{v,k}, \\
    \label{eq:q+}
    q_{k+1} &= A_{qq,k}q_k+A_{q\omega,k}\omega_k+B_{\omega,k}u_k+r_{q,k}, \\
    \label{eq:w+}
    \omega_{k+1} &= A_{\omega\omega,k}\omega_k+B_{\omega,k}u_k+r_{\omega,k},
  \end{align}
\end{subequations}
where \defvar[hide]{$(p_k,v_k,q_k,\omega_k)$}{chaser state at time $kt_c$}
$p_k\definedas p(kt_c)$, $v_k\definedas v(kt_c)$, $q_k\definedas q(kt_c)$ and
$\omega_k\definedas \omega(kt_c)$. The control input is
$u_k=(\Delta t_k^1,\dots,\Delta t_k^M)\in\reals^{16}$, and we shall denote
$u_k^i\definedas\Delta t_k^i$. The expressions for the update matrices
$A_{\cdot\cdot,k}$, $B_{\cdot,k}$ and residual vectors $r_{\cdot,k}$ in
\eqref{eq:dltv} may be obtained directly in our implementation source
code\footnote{In the method \texttt{ApolloCSM.dltv()} in the file
  \texttt{csm.py}.}, a link to which is given in
Section~\ref{subsec:example_parameters}. Given a temporal grid of
\defvar{$N$}{temporal grid density} nodes, the discrete updates~\eqref{eq:dltv}
are imposed at nodes $k=0,\dots,N-1$. This yields a total of $N+1$ states and
$N$ inputs for the discretized dynamics. Physically, each time step in the
discretized dynamics corresponds to a real-time duration of one control
interval.

\subsection{Integer Constraint Handling via State-Triggered Constraints}
\label{subsec:stc}

Even though Problem~\ref{problem:brtop} is finite
dimensional after discretization, it is still a non-convex optimization problem due to the presence of the state-triggered constraints~\eqref{eq:brtop_g}-\eqref{eq:brtop_i}. This type of constraint has recently been
formalized for the method of successive convexification
\cite{Szmuk2018,Szmuk2019a,Szmuk2019b,Reynolds2019a,Reynolds2019b} and is
generally written as:
\begin{equation}
  \label{eq:stc_general}
  g(z(t))<0\Rightarrow c(z(t))\le 0,
\end{equation}
where $z(t)\in\reals^p$ is an arbitrary decision variable (i.e., part of the state, the
control or time), $g(\cdot):\reals^p\to\reals$ is a \textit{trigger function}
and $c(\cdot):\reals^p\to\reals$ is a \textit{constraint function}. Although the
setup can be more general, we shall assume that the constraint function is
convex.

The most direct way of implementing \eqref{eq:stc_general} is by introducing
boolean (integer) variables to imply constraint satisfaction at given times. However, this
leads to a difficult mixed-integer problem with $MN\approx 2000$ boolean
variables. Despite the rapidly advancing state of the art in solution algorithms
and heuristics, mixed-integer programming remains burdened by worst-case
exponential computational complexity \cite{Achterberg2016,Richards2002,Richards2005}. This is
compounded by the iterative nature of successive convexification (see
Section~\ref{subsec:ftr}) which would require solving a mixed-integer program
for several iterations -- typically between five to twenty times.

An alternative approach to incorporating~\eqref{eq:stc_general} is to use continuous variables to formulate a logically equivalent constraint. Recent work \cite{Szmuk2018,Reynolds2019b} has shown
the following equivalence to hold.

\begin{theorem}
  \label{theorem:stc_equivalent}
  The constraint \eqref{eq:stc_general} holds if and only if the following
  constraint holds:
  \begin{equation}
    \label{eq:stc_continuous}
    -\min(g(z(t)),0)c(z(t))\le 0.
  \end{equation}

  \begin{proof}
    Suppose that \eqref{eq:stc_general} holds. If $g(z(t))<0$, then
    $c(z(t))\le 0$ and $\min(g(z(t)),0)c(z(t))=g(z(t))c(z(t))\ge 0$, hence
    \eqref{eq:stc_continuous} holds. If $g(z(t))\ge 0$ then $c(z(t))$ is
    unconstrained. Now suppose that \eqref{eq:stc_continuous} holds. If
    $g(z(t))<0$ then this implies that $g(z(t))c(z(t))\ge 0$, hence it must be
    that $c(z(t))\le 0$, thus \eqref{eq:stc_general} holds. If $g(z(t))\ge 0$
    then $c(z(t))$ is unconstrained.
  \end{proof}
\end{theorem}

The fundamental advantage of \eqref{eq:stc_continuous} is that it is an
equivalent \textit{continuous} variable formulation of \eqref{eq:stc_general}. As a
result, it can be embedded in a non-linear optimization framework without
resorting to integer variables. Using the same reference trajectory $\bar z(\cdot)$ as in the previous section, we can write a convex \textit{successive approximation} of \eqref{eq:stc_continuous}:
\begin{equation}
  \label{eq:stc_approx}
  -\min(g(\bar z(t)),0)c(z(t))\le 0.
\end{equation}

We now approximate the constraints \eqref{eq:brtop_g}-\eqref{eq:brtop_i} in the
form \eqref{eq:stc_approx}:
\begin{subequations}
  \label{eq:stc_traj_approx}
  \begin{align}
    \label{eq:stc_approx_lb}
    &\min(\overline{\Delta} t_k^i-\Delta t_{\min},0)\Delta t_k^i = 0~\text{for
    all}~i=1,\dots,M, \\
    \label{eq:stc_approx_impingement_attitude}
    &-\min(\|\bar p_{k+1}-p_f\|_2-r_a,0)(\cos(\Delta\theta_{\max}/2)-
      e_1^\transp[q_f]_{\otimes}q_k^*)\le 0, \\
    \label{eq:stc_approx_impingement_thrust}
    &\min(\|\bar p_k-p_f\|_2-r_a,0)\Delta t_k^i= 0~\text{for
      all}~i\in\mathcal M.
  \end{align}
\end{subequations}

The constraint set \eqref{eq:stc_traj_approx} is convex in the discrete state and control vectors. Note that \eqref{eq:stc_approx_impingement_attitude} uses a phase lead in the
form of $\bar p_{k+1}$, effectively imposing the impingement constraint on the
attitude one time step in advance. We have found the attitude obtained from
propagating the optimized input trajectory through the non-linear dynamics to
have a slight lag with respect to the optimized attitude trajectory, and this
phase lead helps the actual non-linear dynamics \eqref{eq:vanilla_dynamics} to
satisfy the constraint \eqref{eq:impingement_attitude}. We also mention a
curious property of \eqref{eq:stc_approx_lb} which we call \textit{locking}.

\begin{definition}
  \label{theorem:locking}
  The successive approximation \eqref{eq:stc_approx} \textit{locks} when
  $g(\bar z(t))<0$ implies $g(z(t))<0$. In this situation, $c(z(t))\le 0$ will
  always hold at time $t$ for all subsequent iterations of successive
  convexification.
\end{definition}

\begin{theorem}
  \label{theorem:locking}
  If $g(z)<c(z)$ then the successive approximation \eqref{eq:stc_approx} will
  lock.

  \begin{proof}
    Suppose that $g(\bar z)<0$, then \eqref{eq:stc_approx} implies that
    $c(z)\le 0$. Since $g(z)<c(z)$, this means that $g(z)< 0$. At the next
    iteration, $z$ becomes $\bar z$ and hence $c(z)\le 0$ once again. The
    state-triggered constraint becomes locked.
  \end{proof}
\end{theorem}

We recognize that \eqref{eq:stc_approx_lb} satisfies the conditions of
Theorem~\ref{theorem:locking}, hence it is susceptible to locking. Physically,
this means that once our solution method chooses a pulse width
$\Delta t_k^i<\Delta t_{\min}$, it is guaranteed that $\Delta t_k^i=0$ at all
subsequent successive convexification iterations, and the $i$-th thruster is
guaranteed to remain silent for the $k$-th control interval. A benefit of this
is that this gives a guarantee that when our solution converges, the pulse width
lower bound constraint \eqref{eq:brtop_f} is guaranteed to be satisfied without
any discretization phenomena. A drawback is that by ``committing'' to
$\Delta t_k^i=0$, infeasibility issues may occur. Remarkably, our tests have
shown that this is rarely an issue. Nevertheless, a heuristic recovery method is
incorporated into the algorithm and is discussed in
Section~\ref{subsubsec:reset}.

We note that a large number of other integer constraints can be handled in
identical fashion, such as sensor pointing constraints and further restrictions
on maneuvering capabilities based on keep-out radii. Since our intent is to
demonstrate the capability, we do not explore these cases but only stress that
their inclusion is entirely possible while maintaining computational
tractability. Some other applications of state-triggered constraints can be
found in \cite{Szmuk2018,Szmuk2019a,Szmuk2019b,Reynolds2019a,Reynolds2019b}.

\subsection{Scaling}
\label{subsec:scaling}

Issues associated with floating point precision and the choice of termination
tolerances can make a numerical solver fail when the optimization variables have highly different magnitudes~\cite{Nocedal2006,Gill1981,Liu2017}. The rendezvous problem is susceptible to
this issue since the pulse width control variable can be on the order of 1~ms,
while positions may be on the order to 100~m, which are five orders-of-magnitude different.

Scaling the decision variables attempts to remedy this numerical issue. Possible
strategies include balancing the magnitudes of the dual variables
\cite{Ross2018}, minimizing the condition number of the cost function Hessian at
the solution, or improving the behavior of the cost function first and second
derivatives with respect to machine precision \cite{Gill1981}. We opt for the
standard remedy of applying an affine transformation to non-dimensionalize the
decision variables \cite{Nocedal2006,Gill1981}:
\begin{gather}
  p_k = S_p\hat p_k+s_p,\quad
  v_k=S_v\hat v_k+s_v,\quad
  q_k=\hat q_k,\quad
  \omega_k = S_\omega\hat\omega_k+s_\omega,\quad
  u_k=S_u\hat u_k+s_u.
\end{gather}

Let $({}_{0}{p}_k,{}_{0}{v}_k,{}_{0}{q}_k,{}_{0}{\omega}_k)$ denote the $k$-th state of the first state trajectory, the equations of which are specified in Section~\ref{subsec:ftr}. We define the offset terms as the average value of the corresponding state along the initial guess:
\begin{equation}
  s_p = \frac{1}{N+1}\sum_{i=0}^N{}_{0}{p}_k,\quad
  s_v = \frac{1}{N+1}\sum_{i=0}^N{}_{0}{v}_k,\quad
  s_\omega = \frac{1}{N+1}\sum_{i=0}^N{}_{0}{\omega}_k,\quad
  s_u = \frac{\Delta t_{\max}}{2}\ones{M}.
\end{equation}

The scaling matrices are then computed such that the elements of the
corresponding scaled variable are contained in the $[-1,1]$ interval:
\begin{gather}
  S_p = \diag{\max_{k}|{}_{0}{p}_k-s_p|},~
  S_v = \diag{\max_{k}|{}_{0}{v}_k-s_v|},~
  S_\omega = \diag{\max_{k}|{}_{0}{\omega}_k-s_\omega|},~S_u = \diag{\Delta t_{\max}\ones{M}/2},
\end{gather}
where the $\max|\cdot|$ operation finds, for each dimension, the largest
absolute difference of the trajectory value with the corresponding offset
term. It is useful to define an overall state as
$x_k\definedas (p_k,v_k,q_k,\omega_k)$ and to have a combined affine
transformation:
\begin{equation}
  x_k = S_x\hat x_k+s_x~\text{where}\quad
  S_x=\blkdiag{S_p,S_v,I_4,S_\omega},\quad
  s_x=(s_p,s_v,0_{4\times 1},s_\omega).
\end{equation}

Lastly, the cost \eqref{eq:original_cost} is scaled to maintain a
comparable magnitude to the trust region \eqref{eq:J_tr} and virtual control
\eqref{eq:J_vc} penalties:
\begin{equation}
  J_f = S_J\hat J_f,\quad S_J = \frac{NM\Delta t_{\min}}{4},
\end{equation}
which corresponds to an arbitrary guess that on average a quarter of the
thrusters will be turned on at the minimum pulse width. This choice serves the
purpose of establishing an appropriate relative cost magnitude between $\hat J_f$
and $\hat J_{tr}$ and $\hat J_{vc}$ (the latter two terms are introduced in the next section).

\subsection{Trust Region and Virtual Control Slack Variables}
\label{subsec:slack_variables}

The linearized dynamics \eqref{eq:linearized_differential_form_dynamics} are
a local approximation of the true dynamics \eqref{eq:vanilla_dynamics}
about the reference state and control trajectory. To maintain accuracy,
it is therefore wise to not allow the solution trajectory to deviate too far away from the
reference. Noting that \eqref{eq:vanilla_dynamics} is affine in the control
variable and hence is non-linear only in the state, we impose a \textit{trust
  region} constraint which keeps the state trajectory close to the reference:
\defvar[hide]{$\hat\eta_k$}{quadratic trust region size at the $k$-th time step}
\begin{equation}
  \label{eq:trust_region}
  \|\hat x_k-S_x^{-1}(\bar x_k-s_x)\|_2^2\le\hat\eta_k,\quad k=1,\dots,N-1,
\end{equation}
where the trust region is not imposed at the boundary nodes where the state is
fixed. Note that by considering the deviation of scaled states, the trust region
is automatically a scaled quantity. The size of the trust region is penalized
with a (large) weight $w_{tr}\in\reals_{\ge 0}$ via an additional cost, such
that the optimization procedure is encouraged to keep the trust region small:
\defvar[hide]{$w_{tr}$}{trust region penalty weight}
\begin{equation}
  \label{eq:J_tr}
  \hat J_{tr} = w_{tr}\sum_{i=1}^{N-1}\hat\eta_k.
\end{equation}

Another issue that arises when moving away from the reference trajectory while
using linearized dynamics is that of artificial infeasibility
\cite{Nocedal2006,Mao2016,Mao2017,Mao2019}. This occurs when the linearized
dynamics (which are not the \textit{true} dynamics) cannot be satisfied while
respecting all the other constraints. The standard remedy to artificial
infeasibility is to introduce a \textit{virtual control} term whose usage
is highly penalized:
\begin{subequations}
  \label{eq:dltv_with_vc}
  \begin{align}
    \label{eq:p+_vc}
    p_{k+1} &= p_k+A_{pv,k}v_k+A_{pq,k}q_k+A_{p\omega,k}\omega_k+B_{p,k}u_k+r_{p,k}+\ell_k^p, \\
    \label{eq:v+_vc}
    v_{k+1} &= v_k+A_{vq,k}q_k+A_{v\omega,k}\omega_k+B_{v,k}u_k+r_{v,k}+\ell_k^v, \\
    \label{eq:q+_vc}
    q_{k+1} &= A_{qq,k}q_k+A_{q\omega,k}\omega_k+B_{\omega,k}u_k+r_{q,k}+\ell_k^q, \\
    \label{eq:w+_vc}
    \omega_{k+1} &= A_{\omega\omega,k}\omega_k+B_{\omega,k}u_k+r_{\omega,k}+\ell_k^\omega,
  \end{align}
\end{subequations}
where $\ell_k\definedas(\ell_k^p,\ell_k^v,\ell_k^q,\ell_k^\omega)\in\reals^{13}$
is the virtual control. We scale the virtual control similarly to the state but
without centering, such that the ideal scaled virtual control is zero:
\begin{equation}
  \label{eq:vc_scaling}
  \ell_k = S_{\ell}\hat \ell_k,\quad S_{\ell} = \diag{\max_{k}|{}_{0}{x_k}|}.
\end{equation}

Use of virtual control is penalized with an additional cost term that has a
(large) weight $w_{vc}\in\reals_{\ge 0}$: \defvar[hide]{$w_{vc}$}{virtual
  control penalty weight}
\begin{equation}
  \label{eq:J_vc}
  \hat J_{vc} = w_{vc}\sum_{i=0}^{N-1}\|\hat{\ell}_k\|_1.
\end{equation}


\subsection{Local Rendezvous Guidance Problem}
\label{subsec:cvx}

We are now in a position to formulate a finite-dimensional convex approximation
of Problem~\ref{problem:brtop}, which we call the \defabbrev{LRGP}{local
  rendezvous trajectory problem}. Solving LRGP to global optimality is the
core step in the successive convexification solution method that is described in the
next section. Note that LRGP relies on a previous solution to itself from a
previous iteration of successive convexification, whose values are demarcated
with a bar. The state, input and virtual control variables are understood to be
substituted with their scaled expressions from Section~\ref{subsec:scaling}.

\begin{poptimization}[Local Rendezvous Guidance
  Problem]{center}{lrtop}{\Delta t_{k=0,\dots,N-1}^{i=1,\dots,M}}{\hat J_f+\hat J_{vc}+\hat J_{tr}}
  p_{k+1} = p_k+A_{pv,k}v_k+A_{pq,k}q_k+A_{p\omega,k}\omega_k+B_{p,k}u_k+r_{p,k}+\ell_k^p,~&&k=0,\dots,N-1, \\
  v_{k+1} = v_k+A_{vq,k}q_k+A_{v\omega,k}\omega_k+B_{v,k}u_k+r_{v,k}+\ell_k^v,~&&k=0,\dots,N-1, \\
  q_{k+1} = A_{qq,k}q_k+A_{q\omega,k}\omega_k+B_{\omega,k}u_k+r_{q,k}+\ell_k^q,~&&k=0,\dots,N-1, \\
  \omega_{k+1} =
  A_{\omega\omega,k}\omega_k+B_{\omega,k}u_k+r_{\omega,k}+\ell_k^\omega,~&&k=0,\dots,N-1,
  \\
  0\le\Delta t_k^i\le\Delta t_{\max}~\text{for
    all}~i=1,\dots,M,~&&k=0,\dots,N-1, \\
  \min(\overline{\Delta} t_k^i-\Delta t_{\min},0)\Delta t_k^i = 0~\text{for
    all}~i=1,\dots,M,~&&k=0,\dots,N-1, \\
  -\min(\|\bar p_{k+1}-p_f\|_2-r_a,0)(\cos(\Delta\theta_{\max}/2)-
  e_1^\transp[q_f]_{\otimes}q_k^*)\le 0,~&&k=1,\dots,N-1, \\
  \min(\|\bar p_k-p_f\|_2-r_a,0)\Delta t_k^i= 0~\text{for all}~i\in\mathcal M,
  ~&&k=0,\dots,N-1,
  \\
  \|p_k-p_l\|_2\cos(\gamma)\le (p_k-p_l)^\transp\hat e_d,~&&k=1,\dots,N-1, \\
  \|\hat x_k-S_x^{-1}(\bar x_k-s_x)\|_2^2\le\hat\eta_k~&&k=1,\dots,N-1, \\
  p_0=p_0,~v_0=v_0,~q_0=q_0,~\omega_0=\omega_0, \\
  p_N=p_f,~v_N=v_f,~q_N=q_f,~\omega_N=\omega_f.
\end{poptimization}

\subsection{Successive Convexification}
\label{subsec:ftr}

\begin{figure}
  \centering
  \includegraphics[width=1\textwidth]{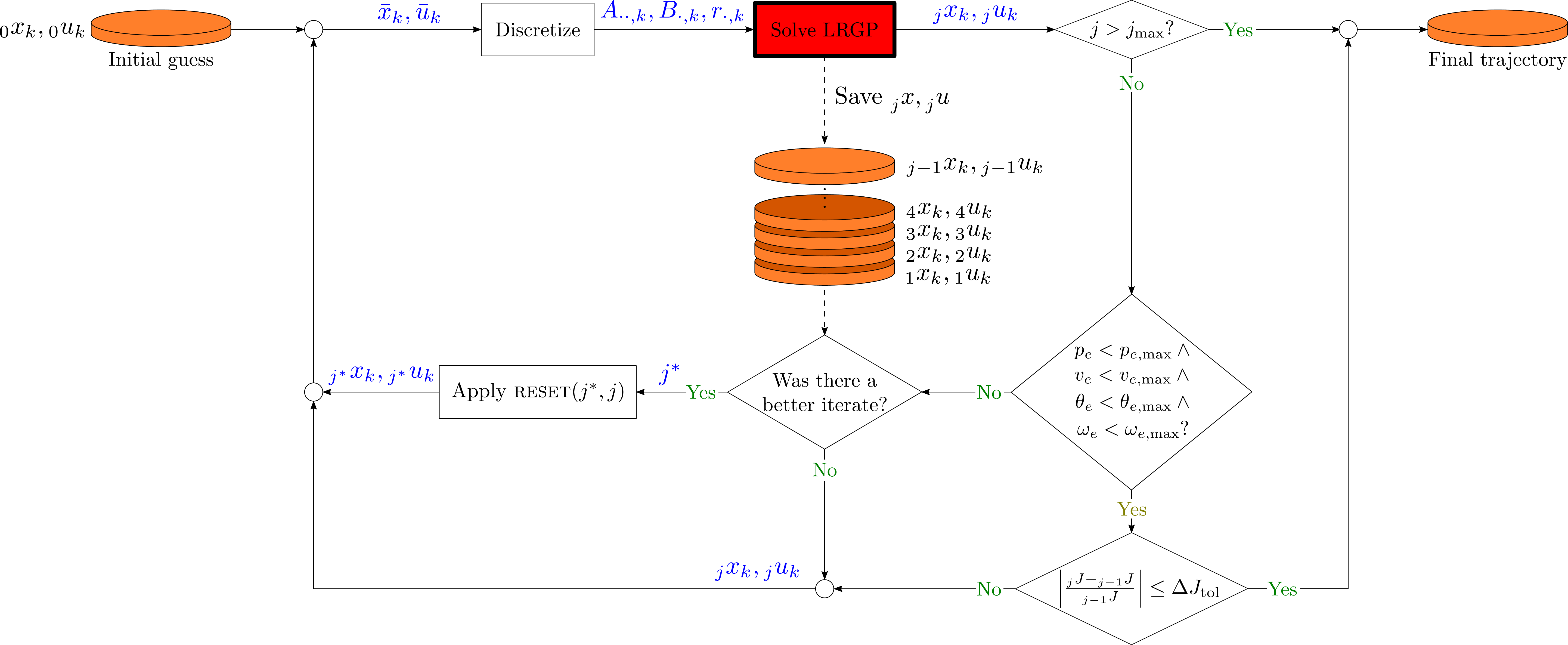}
  \caption{Block diagram illustration of the successive convexification
    procedure. Algorithm~\ref{alg:ftr} provides the corresponding pseudocode
    implementation. The core step is the solution of
    Problem~\ref{problem:lrtop}, which is a second-order cone program.}
  \label{fig:successive_convexification_block_diagram}
\end{figure}

We solve Problem~\ref{problem:brtop} using a successive
convexification methodology. The solutions that are obtained with this approach have the
following properties \cite{Mao2019,Reynolds2019b}:
\begin{itemize}
\item The final trajectory exactly satisfies the non-linear dynamics
  \eqref{eq:vanilla_dynamics};
\item The input constraints \eqref{eq:brtop_f} and \eqref{eq:brtop_g} are
  exactly satisfied;
\item The state constraints \eqref{eq:brtop_h}-\eqref{eq:brtop_j} are satisfied
  at each discrete-time grid node $k=0,\dots,N$.\footnote{By making the temporal
    grid refined enough, or by incorporating safety margins of a large enough
    size, the state constraints \eqref{eq:brtop_h}-\eqref{eq:brtop_j} can be
    guaranteeably satisfied for all times $t\in [0,t_f]$. A similar procedure is
    done in \eqref{eq:stc_approx_impingement_attitude}, where phase lead helps
    to achieve continuous-time constraint satisfaction.}
\end{itemize}

It remains an open question whether the converged solution of the successive
convexification flavor presented in this work, called the \textit{penalized trust
  region} method, is a local optimum of the original
Problem~\ref{problem:brtop}. What is known is that, because there is incentive to
minimize $\hat J_f$, the converged rendezvous trajectory is fuel efficient. In particular, Section~\ref{sec:example} shows that we are able to generate trajectories that are far more fuel efficient that those executed during an actual Apollo rendezvous maneuver.

The successive convexification method is illustrated in
Figure~\ref{fig:successive_convexification_block_diagram} and works by
iteratively solving Problem~\ref{problem:lrtop}. Henceforth, we adopt the
convention that a variable like ${}_{j}p_k$ represents the solution obtained at
the $j$-th successive convexification iteration and the $k$-th discrete-time
index. Like most non-linear optimization routines, the algorithm requires an
initial trajectory guess. However, successive convexification is known to be
very flexible in the quality of the initial guess that it is
given~\cite{Szmuk2018,Reynolds2019b,Bonalli2019}. As a result, we propose a
constant-velocity linearly interpolated initial trajectory between the initial
and final states:
\begin{equation}
  \label{eq:initial_state_trajectory_guess}
  {}_{0}{p_k} = p_0+\frac{k}{N}p_f,\quad
  {}_{0}{v_k} = \frac{p_f-p_0}{t_f},\quad
  {}_{0}{q_k} = q_0\otimes
  \begin{bmatrix}
    \cos\left(\frac{k\theta_{0f}}{2N}\right) \\
    u_{0f}\sin\left(\frac{k\theta_{0f}}{2N}\right)
  \end{bmatrix},\quad
  {}_{0}{\omega_k} = \frac{\theta_{0f}}{N}u_{0f},
\end{equation}
where
$q_0^*\otimes q_f=(\cos(\theta_{0f}/2),u_{0f}\sin(\theta_{0f}/N))\in\quaternion$
is the error quaternion between the initial and final quaternions such that
${}_{0}{q_k}$ corresponds to a spherical linear interpolation
\cite{Sola2017}. We also use a constant minimum pulse width initial input
trajectory:
\begin{equation}
  \label{eq:initial_input_trajectory_guess}
  {}_{0}{u_k} = \Delta t_{\min},
\end{equation}
which ensures that the constraint \eqref{eq:brtop_g} does not lock from the
outset of the solution process, as per the discussion at the end of
Section~\ref{subsec:stc}. Note that neither does
\eqref{eq:initial_state_trajectory_guess} satisfy the dynamics
\eqref{eq:vanilla_dynamics} nor does \eqref{eq:initial_input_trajectory_guess}
generate this trajectory. The fact that this is an acceptable initial guess for
the successive convexification method is testimony to the method's
aforementioned insensitivity to the initialization quality.

With an initialization available,
Figure~\ref{fig:successive_convexification_block_diagram} shows how each
iteration of successive convexification first linearizes and discretizes the
dynamics \eqref{eq:vanilla_dynamics} about a previously obtained state and input
trajectory, and then solves Problem~\ref{problem:lrtop}. The discretization
method was described in Section~\ref{subsec:propagation} and requires
\textit{continuous-time} state and input trajectories. However, the optimization
problem only outputs trajectories that are defined at discrete-time nodes. To
reconcile these facts, we propagate the discrete-time input through the
non-linear dynamics \eqref{eq:vanilla_dynamics} by using the discrete-time input
trajectory that is converted to continuous-time via the thruster model
\eqref{eq:thrust_model}. Letting
$\bar x_k=(\bar p_k,\bar v_k,\bar q_k,\bar\omega_k)$ denote the discrete-time
reference state trajectory at the $k$-th time step, we have for
$t\in [kt_c,(k+1)t_c]$:
\begin{subequations}
  \label{eq:propagation_with_resetting}
  \begin{align}
    \bar p(t) &= \bar p_k+\int_{kt_c}^{t}\bar v(\tau)\dd\tau, \\
    \bar v(t) &= \bar v_k+\frac{1}{m}\sum_{i=1}^M\int_{kt_c}^t\bar q(\tau)\otimes
                \bar f_i(\tau)\otimes\bar q(\tau)^*\dd\tau, \\
    \bar q(t) &= \bar q_k+\frac{1}{2}\int_{kt_c}^{t}\bar q(\tau)\otimes\bar\omega(\tau)\dd\tau, \\
    \bar\omega(t) &= \bar\omega_k+J^{-1}\int_{kt_c}^t\biggl[
                    \sum_{i=1}^M r_i\times\bar f_i(\tau)-
                    \bar\omega(\tau)\times (J\bar\omega(\tau))\biggr]\dd\tau,
  \end{align}
\end{subequations}

where, following from the thruster model \eqref{eq:thrust_model}:
\begin{equation}
  \bar f_i(t) = 
  \begin{cases}
    \hat f_i & \text{if }t\in [k t_c,k t_c+\bar{\Delta} t_k^i], \\
    0 & \text{otherwise}.
  \end{cases}  
\end{equation}

Two aspects about \eqref{eq:propagation_with_resetting} deserve
attention. First, the pulsed on-off nature of the thrust signal makes
\eqref{eq:propagation_with_resetting} a stiff integration problem. However,
since we know that the falling edge of the $i$-th thrust signal occurs at
$\bar{\Delta} t_k^i$, we can take this explicitly into account such that a stiff
numerical integrator is not required. For example, if we assume (without loss of
generality) that $\bar{\Delta} t_k^i$ is the $i$-th shortest pulse width, then
we can integrate over $[kt_c,kt_c+\bar{\Delta} t_k^1]$, then over
$[kt_c+\bar{\Delta} t_k^1,kt_c+\bar{\Delta} t_k^2]$, and proceed accordingly
until $[kt_c+\bar{\Delta} t_k^{M-1},kt_c+\bar{\Delta} t_k^M]$ and, finally, over
$[kt_c+\bar{\Delta} t_k^M,(k+1)t_c]$. In implementation, we use the adaptive
stepsize Dormand-Price integration method, which is an explicit Runge-Kutta
method where integration error is controlled by assuming the accuracy of a
fourth-order solution and where steps are taken using the fifth-order solution
via local extrapolation \cite{Dormand1980}.

\begin{figure}
  \centering
  \includegraphics[width=0.5\textwidth]{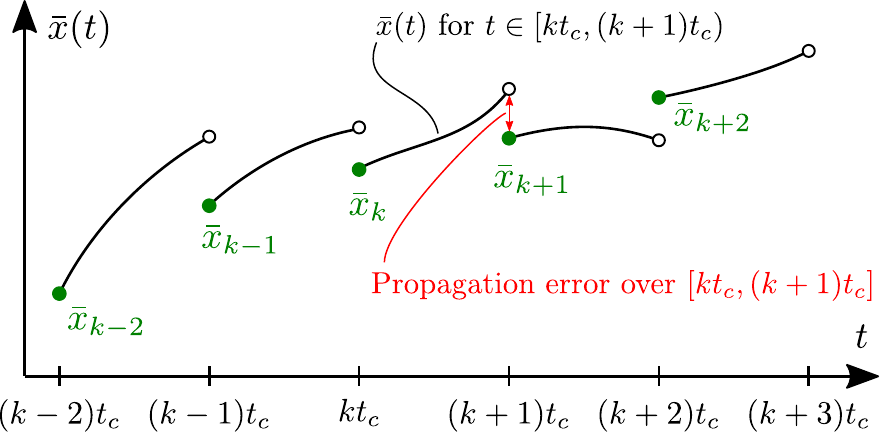}
  \caption{Illustration of the reference state trajectory obtained via
    \eqref{eq:propagation_with_resetting}. Reinitialization at the beginning of
    each control interval creates discontinuities which we call the propagation
    error. Reinitialization improves successive convexification convergence.}
  \label{fig:propagation}
\end{figure}

The second subtlety of \eqref{eq:propagation_with_resetting} is that the
integration is \textit{reinitialized} over each control interval to the
discrete-time reference state $\bar x_k$. This leads to a discontinuous
reference state trajectory as illustrated in Figure~\ref{fig:propagation}. At
the end of each control interval, the discrepancy between the non-linear
propagation and the state variable returned by Problem~\ref{problem:lrtop} is
termed the \textit{propagation error}.

\begin{definition}
  \label{definition:propagation_error}
  We call the \textit{propagation error} over a time interval $[k_1t_c,k_2t_c]$,
  $k_2>k_1$, the difference $\bar x_{k_2}-\bar x(k_2t_c)$, where:
  \begin{enumerate}[label=(\arabic*)]
  \item $\bar x_{k_2}$ is the state returned by Problem~\ref{problem:lrtop} for time
    step $k_2\in\integers_{> 0}$;
  \item $\bar x(k_2t_c)$ is the state resulting from integrating the non-linear
    dynamics \eqref{eq:vanilla_dynamics} initialized at state $\bar x_{k_1}$ and
    using the input trajectory returned by Problem~\ref{problem:lrtop} over the
    $[k_1t_c,k_2t_c]$ interval, as is done in \eqref{eq:propagation_with_resetting}
    over the $[kt_c,(k+1)t_c]$ interval.
  \end{enumerate}
\end{definition}

We have found that by reinitializing the integration process over each control
interval, the convergence of successive convexification is improved since
the propagation error is kept small. The benefit is analogous to how multiple-shooting
is more stable than single-shooting solution methods
\cite{Szmuk2018,Malyuta2019}.

Complementary to Figure~\ref{fig:successive_convexification_block_diagram},
Algorithm~\ref{alg:ftr} formalizes the successive convexification method for
rendezvous trajectory generation. We draw attention to several special
features. First, a maximum iteration count \defvar{$j_{\max}$}{maximum iteration
  count} bounds the maximum runtime, where the index $j$ is reserved for the
successive convexification algorithm's iteration number. Second, the algorithm
will terminate if the cost \eqref{eq:original_cost} does not change by more than
a relative tolerance \defvar{$\Delta J_{\text{tol}}$}{relative cost change
  tolerance}. Lastly, we introduce a novel \textit{resetting} feature which
partially handles feasibility issues associated with locking of the
implementation \eqref{eq:lrtop_g} of the input constraint
\eqref{eq:brtop_g}. The next subsection explains this procedure in detail.

\begin{algorithm}
  \centering
  \begin{algorithmic}[1]
    \State
    ${}_{0}{x_k},{}_{0}{u_k}\gets\text{Initial guess
      \eqref{eq:initial_state_trajectory_guess} and
      \eqref{eq:initial_input_trajectory_guess}}$
    \State $\bar x_k,\bar u_k\gets {}_{0}{x},{}_{0}{u}$
    \Comment{Initialize reference trajectories to the initial guess}
    \For{$j=1,\dots,j_{\max}+1$}
    \If{$j>j_{\max}$}
    \State \textbf{break} \Comment{Reached maximum iterations}
    \Else
    \State \begin{varwidth}[t]{\linewidth}
      $A_{pv,k},A_{pq,k},A_{p\omega,k},A_{vq,k},A_{v\omega,k},A_{qq,k},A_{q\omega,k},
      A_{\omega\omega,k},B_{p,k},B_{v,k},B_{q,k},B_{\omega,k},r_{p,k},r_{v,k},r_{q,k},
      r_{\omega,k}\gets$ Linearize\par
      \hskip\algorithmicindent
      about the $\bar x_k$, $\bar u_k$
      reference and discretize
    \end{varwidth}
    \State ${}_{j}{J},{}_{j}{x_k},{}_{j}{u_k}\gets\text{Solve
      Problem~\ref{problem:lrtop}}$ \Comment{Core step (call to convex
      optimizer)}
    \label{alg:ftr:line:solve}
    \State Save ${}_{j}{J},{}_{j}{x_k},{}_{j}{u_k}$ into a list of iterates
    \State
    $\bar p(t_f),\bar v(t_f),\bar q(t_f),\bar\omega(t_f)\gets\text{Compute
      \eqref{eq:propagation_with_resetting} using $k=0$ and $t=t_f$}$
    \label{alg:ftr:line:accuracy_check_start}
    \Comment{Propagation error over $[0,t_f]$}
    \State $p_e\gets \|{}_{j}{p_N}-\bar p(t_f)\|_{\infty}$
    \State $v_e\gets \|{}_{j}{v_N}-\bar v(t_f)\|_{\infty}$
    \State
    $\omega_e\gets \frac{180}{\pi}\|{}_{j}{\omega_N}-\bar \omega(t_f)\|_{\infty}$
    \State
    $\theta_e\gets\frac{360}{\pi}\arccos(e_1^\transp({}_{j}{q_N^*}\otimes\bar q(t_f)))$
    \If{${}_{j}{p_e}<p_{e,\max}\wedge {}_{j}{v_e}<v_{e,\max}\wedge
      {}_{j}{\theta_e}<\theta_{e,\max}\wedge
      {}_{j}{\omega_e}<\omega_{e,\max}$}
    \label{alg:ftr:line:accuracy_check_end}
    \If{$\left|\frac{{}_{j}{J}-{}_{j-1}{J}}
        {{}_{j-1}{J}}\right|\le\Delta J_{\text{tol}}$}
    \State \textbf{break} \Comment{Converged}
    \EndIf
    \State $\bar x_k,\bar u_k\gets {}_{j}{x_k},{}_{j}{u_k}$
    \Else
    \State $j^*\gets\max_{l=1,\dots,j-1}\{l\mid
    {}_{l}{p_e}<p_{e,\max}\wedge {}_{l}{v_e}<v_{e,\max}\wedge
    {}_{l}{\theta_e}<\theta_{e,\max}\wedge
    {}_{l}{\omega_e}<\omega_{e,\max}\wedge\text{\eqref{eq:pulse_width_constraint} holds}\}$
    \label{alg:ftr:line:find_j_star_start}
    \If{no $j^*$ found}
    \State $j^*\gets\max_{l=1,\dots,j-1}\{l\mid
    {}_{l}{p_e}<p_{e,\max}\wedge {}_{l}{v_e}<v_{e,\max}\wedge
    {}_{l}{\theta_e}<\theta_{e,\max}\wedge
    {}_{l}{\omega_e}<\omega_{e,\max}\}$
    \label{alg:ftr:line:find_j_star_end}
    \EndIf
    \If{no $j^*$ found}
    \State $\bar x_k,\bar u_k\gets {}_{j}{x_k},{}_{j}{u_k}$
    \Comment{Continue with what we have}
    \Else 
    \State $\bar x_k,\bar u_k\gets\text{\Call{reset}{$j^*$,$j$}}$
    \EndIf
    \EndIf
    \EndIf
    \EndFor
    \Statex
    \Function{reset}{$j^*,j$}
    \label{alg:ftr:line:reset_start}
    \State $x,u\gets{}_{j^*}{x_k},{}_{j^*}{u_k}$
    \For{$k=0,\dots,N-1$}
    \For{$i=1,\dots,M$}
    \For{$l=j^*+1,\dots,j-1$}
    \If{${}_{l}{\Delta t_k^i}\in (0,\Delta t_{\min})$}
    \State Add the constraint $\Delta t_k^i \ge \Delta t_{\min}$ to
    Problem~\ref{problem:lrtop}
    \State $u_k^i\gets\Delta t_{\min}$ \Comment{Reset ${}_{j^*}\Delta t_k^i$ to
      $\Delta t_{\min}$}
    \EndIf
    \EndFor
    \EndFor
    \EndFor
    \State \textbf{return} $x,u$
    \label{alg:ftr:line:reset_end}
    \EndFunction
  \end{algorithmic}
  \caption{Successive convexification algorithm for rendezvous trajectory
    generation. This finds a feasible solution to Problem~\ref{problem:brtop}
    which in practice also has a low enough cost to be useful for engineering purposes.}
  \label{alg:ftr}
\end{algorithm}

\subsubsection{Solution Resetting}
\label{subsubsec:reset}

We define a solution obtained during iteration $j$ of successive convexification
as \textit{feasible} if the propagation error over $[0,t_f]$ is smaller than a
given tolerance. In particular, a user specifies a position accuracy
$p_{e,\max}$, a velocity accuracy $v_{e,\max}$, an angular error accuracy $\theta_{e,\max}$, and an angular velocity accuracy $\omega_{e,\max}$. The propagation error is compared against these accuracy values on lines~\ref{alg:ftr:line:accuracy_check_start}-\ref{alg:ftr:line:accuracy_check_end} of Algorithm~\ref{alg:ftr}. If a particular solution does not satisfy this
accuracy, then we perform a \textit{solution reset}. This becomes particularly
useful when locking of \eqref{eq:brtop_g} leads to an infeasibility, thereby
prompting virtual control usage and yielding a large propagation error. In
effect, the solution at iteration $j$ becomes \textit{worse} than an earlier
iterate. When such a degradation in solution accuracy occurs, we first search
for a ``better'' earlier iterate $j^*$ using the following logic:
\begin{itemize}
\item We first try to find an iterate whose solution is accurate to within the
  required tolerance \textit{and} whose input trajectory satisfies
  \eqref{eq:pulse_width_constraint};
\item If no iterate satisfies \eqref{eq:pulse_width_constraint}, we simply try
  to find an iterate whose solution is accurate enough.
\end{itemize}

If such a $j^*$ is found, a reset operation is performed by the \textsc{reset}
function on
lines~\ref{alg:ftr:line:reset_start}-\ref{alg:ftr:line:reset_end}. The result is
to lower-bound the pulse width to $\Delta t_{\min}$ for all $\Delta t_k^i$
values which violate \eqref{eq:pulse_width_constraint} at iterations after $j^*$
and before $j$. The motivation for doing this can be seen in two steps:
\begin{enumerate}
\item All such $\Delta t_k^i$'s become locked;
\item At least one of the $\Delta t_k^i$'s that locked now needs to become
  positive in order to avoid using the virtual control. Since locking by
  definition prevents this, virtual control is used and solution accuracy
  degrades;
\item By forbidding all such $\Delta t_k^i$'s from being lower than
  $\Delta t_{\min}$ we (conservatively) prevent locking in a subset of
  $\Delta t_k^i$ values whose locking has lead to virtual control usage.
\end{enumerate}

From the above reasoning, it is clear why on
line~\ref{alg:ftr:line:find_j_star_start} we first try to search for a $j^*$
iterate where \eqref{eq:pulse_width_constraint} holds. Roughly speaking, the
solution for such an iterate is \textit{already} a feasible solution and hence
the \textsc{reset} function can only exploit it to obtain a feasible solution
that is, in the worst case, the same as the solution at $j^*$. This is
formalized in the following theorem.

\begin{theorem}
  \label{label:reset_worst_case}
  Suppose that an iterate $j^*$ is found on
  line~\ref{alg:ftr:line:find_j_star_start}. Assume that $w_{vc}$ is
  sufficiently large and that infeasibility due to re-linearization of the
  system about ${}_{j^*}x_k$ (i.e. artificial infeasibility as discussed in
  Section~\ref{subsec:slack_variables}) does not occur. By using the
  \textsc{reset} scheme, after enough iterations (assume $j_{\max}=\infty$) it
  is guaranteed that the solution returned by successive convexification
  satisfies \eqref{eq:pulse_width_constraint} and the propagation error
  tolerance.

  \begin{proof}
    The solution at $j^*$ already satisfies \eqref{eq:pulse_width_constraint}
    and the propagation error tolerance. Since no ${}_{j^*}\Delta t_k^i$
    violates \eqref{eq:pulse_width_constraint}, we have either
    ${}_{j^*}\Delta t_k^i\ge\Delta t_{\min}$ or ${}_{j^*}\Delta t_k^i=0$. The
    effect of repeated application of \textsc{reset} is that eventually all
    non-zero $\Delta t_k^i$ inputs are lower-bounded by $\Delta t_{\min}$. At
    this point, locking can no longer occur and the only degrees of freedom that
    remain in the optimization are to adjust the non-zero pulse widths in the
    $[\Delta t_{\min},\Delta t_{\max}]$ band. For a large enough $w_{vc}$, the
    penalty of increasing virtual control will outweigh any benefits of
    deviating $\Delta t_k^i$ from ${}_{j^*}\Delta t_k^i$. As a result,
    successive convexification will, in the worst case, return a solution that
    also satisfies \eqref{eq:pulse_width_constraint} as well as the propagation
    error tolerance.
  \end{proof}
\end{theorem}

Note that Theorem~\ref{label:reset_worst_case} relies on two important
qualifiers: a ``large enough'' $w_{vc}$ and the absence of artificial
infeasibility. The former is a loose statement, however we have found in our
numerical studies in Section~\ref{sec:example} that an appropriate $w_{vc}$ is
very easy to find as some (large) power of 10. The latter qualifier is a
fictitious one made for the purpose of the proof, which motivates the use of the
\textsc{reset} \textit{heuristic}. Indeed, because the solution obtained at
$j^*$ is for a linearized system, it is still possible that infeasibility after
$j^*$ occurs not due to locking but due to the effect re-linearizing. This is a
side-effect of taking too large a step during the iteration process. The remedy
is to increase $w_{tr}$ such that the trust region is made smaller. In our
experience, the \textsc{reset} function together with a large enough $w_{tr}$
setting leads to a quite reliable solution method. Like $w_{vc}$, a ``large
enough'' $w_{tr}$ is easy of find as a power of 10.

\section{Apollo Transposition and Docking Example}
\label{sec:example}

In this section we apply our solution method to a real docking maneuver between the Apollo \defabbrev{CSM}{Command and Service Module} and the \defabbrev{LM}{Lunar
  Module}. Section~\ref{subsec:example_parameters} defines the problem
parameters and Section~\ref{subsec:results} discusses the resulting rendezvous
trajectory and its properties. The main takeaway is that we are able to generate
in several seconds a trajectory that is up to $90~\%$ more fuel optimal than the
Apollo design target.

\subsection{Problem Parameters}
\label{subsec:example_parameters}

\begin{figure}
  \centering
\begin{tikzpicture}[scale=0.7]
\tikzmath{
	\cx1=-3; \cy=0; \wdth=5; \hght=4;
	\cx2=\cx1+\wdth+0.5;
	\cx3=\cx2+\wdth+0.5;
	\cx4=\cx3+\wdth+0.5;
	\svx=-5; \svy=0;
	\svxx=0.5; \svyy=0;
	\svxxx=6; \svyyy=0;
	\svxxxx=11.5; \svyyyy=0;
	\txx=cos(-20); \tyy=sin(-20);
}

\pane{\cx1}{\cy}{0}{\wdth}{\hght}{0.5}{color=black,fill=beige!60}
\pane{\cx2}{\cy}{0}{\wdth}{\hght}{0.5}{color=black,fill=beige!60}
\pane{\cx3}{\cy}{0}{\wdth}{\hght}{0.5}{color=black,fill=beige!60}
\pane{\cx4}{\cy}{0}{\wdth}{\hght}{0.5}{color=black,fill=beige!60}

\begin{scope}[shift={(\svx,\svy)},rotate=-20]
\coordinate (A) at (-0.53,0);
\coordinate (B) at ($(A)+(0.75,0)$);
\coordinate (C) at ($(B)+(0,1.25)$);
\draw[thick,shading=axis,bottom color=darkgrey!50!white,top color=darkgrey]  (A) -- ++(0.75,0) to[bend left] ++(0,1.25) -- ++(-1.21,0);
\draw[shading=axis,bottom color=darkgrey!80!white,top color=darkgrey] (C) to[out=20,in=190] ++(0.8,0.2) to[bend right] ++(0,-0.1) to[out=180,in=20] ($(C)+(0,-0.05)$) -- cycle;
\draw[shading=axis,bottom color=darkgrey!80!white,top color=darkgrey] ($(B)+(0,0.05)$) to[out=-20,in=175] ++(0.8,-0.2) to[bend right] ++(0,-0.1) to[out=170,in=-20] (B) -- cycle;
\draw[shading=ball,inner color=darkgrey,outer color=darkgrey!70!white] (B) -- ++(0.6,0.4) to[bend left] ++(0,0.4) -- (C) to[bend right] cycle;
\end{scope}
\apolloCSM{-3.5}{0.1}{-110}{0}{0.25}
\draw[ultra thick,->] (-1.9,-0.485) -- ++(\txx,\tyy);

\begin{scope}[shift={(\svxx,\svyy)},rotate=-20]
\coordinate (A) at (-0.53,0);
\coordinate (B) at ($(A)+(0.75,0)$);
\coordinate (C) at ($(B)+(0,1.25)$);
\draw[thick,shading=axis,bottom color=darkgrey!50!white,top color=darkgrey]  (A) -- ++(0.75,0) to[bend left] ++(0,1.25) -- ++(-1.21,0);
\draw[shading=axis,bottom color=darkgrey!80!white,top color=darkgrey] (C) to[out=20,in=190] ++(0.8,0.2) to[bend right] ++(0,-0.1) to[out=180,in=20] ($(C)+(0,-0.05)$) -- cycle;
\draw[shading=axis,bottom color=darkgrey!80!white,top color=darkgrey] ($(B)+(0,0.05)$) to[out=-20,in=175] ++(0.8,-0.2) to[bend right] ++(0,-0.1) to[out=170,in=-20] (B) -- cycle;
\draw[shading=ball,inner color=darkgrey,outer color=darkgrey!70!white] (B) -- ++(0.6,0.4) to[bend left] ++(0,0.4) -- (C) to[bend right] cycle;
\end{scope}
\apolloCSM{3}{-1}{0}{0}{0.25}
\draw[ultra thick,->] (3.5,-1.25) arc (-90:90:0.8); 

\begin{scope}[shift={(\svxxx,\svyyy)},rotate=-20]
\coordinate (A) at (-0.53,0);
\coordinate (B) at ($(A)+(0.75,0)$);
\coordinate (C) at ($(B)+(0,1.25)$);
\draw[thick,shading=axis,bottom color=darkgrey!50!white,top color=darkgrey]  (A) -- ++(0.75,0) to[bend left] ++(0,1.25) -- ++(-1.21,0);
\draw[shading=axis,bottom color=darkgrey!80!white,top color=darkgrey] (C) to[out=20,in=190] ++(0.8,0.2) to[bend right] ++(0,-0.1) to[out=180,in=20] ($(C)+(0,-0.05)$) -- cycle;
\draw[shading=axis,bottom color=darkgrey!80!white,top color=darkgrey] ($(B)+(0,0.05)$) to[out=-20,in=175] ++(0.8,-0.2) to[bend right] ++(0,-0.1) to[out=170,in=-20] (B) -- cycle;
\draw[shading=ball,inner color=darkgrey,outer color=darkgrey!70!white] (B) -- ++(0.6,0.4) to[bend left] ++(0,0.4) -- (C) to[bend right] cycle;
\end{scope}
\apolloCSM{8.75}{-0.375}{70}{0}{0.25}
\draw[ultra thick,<-] (9.2,-0.525) -- ++(\txx,\tyy);

\begin{scope}[shift={(\svxxxx,\svyyyy)},rotate=-20]
\coordinate (A) at (-0.53,0);
\coordinate (B) at ($(A)+(0.75,0)$);
\coordinate (C) at ($(B)+(0,1.25)$);
\draw[thick,shading=axis,bottom color=darkgrey!50!white,top color=darkgrey]  (A) -- ++(0.75,0) to[bend left] ++(0,1.25) -- ++(-1.21,0);
\draw[shading=axis,bottom color=darkgrey!80!white,top color=darkgrey] (C) to[out=20,in=190] ++(0.8,0.2) to[bend right] ++(0,-0.1) to[out=180,in=20] ($(C)+(0,-0.05)$) -- cycle;
\draw[shading=axis,bottom color=darkgrey!80!white,top color=darkgrey] ($(B)+(0,0.05)$) to[out=-20,in=175] ++(0.8,-0.2) to[bend right] ++(0,-0.1) to[out=170,in=-20] (B) -- cycle;
\draw[shading=ball,inner color=darkgrey,outer color=darkgrey!70!white] (B) -- ++(0.6,0.4) to[bend left] ++(0,0.4) -- (C) to[bend right] cycle;
\end{scope}
\apolloCSM{15.04}{-0.64}{70}{0}{0.25}
\apollodocked{13}{0.1}{0.2}{-110}{0}{0}

\node[align=center, text width=5cm] at (-3,-2.55) {\small CSM separates from \\ the S-IVB};
\node[align=center, text width=5cm] at (2.5,-2.55) {\small CSM executes a \\ $180^{\circ}$ rotation};
\node[align=center, text width=5cm] at (8,-2.55) {\small CSM closes and aligns with berthed LM};
\node[align=center, text width=5cm] at (13.5,-2.6) {\small CSM/LM dock and \\ separate from the S-IVB};
\draw[thick,|<->|] (0,-3.5) -- ++(10.5,0) node[anchor=north,pos=0.5] {\small Successive convexification is applied to these phases};

\end{tikzpicture}
  \caption{Illustration of the Apollo CSM transposition and docking maneuver
    with the LM housed inside the S-IVB \cite[Figure~2-11]{apollo_mass}. We
    design an open-loop reference trajectory for the middle two phases.}
  \label{fig:transposition_docking_illustration}
\end{figure}
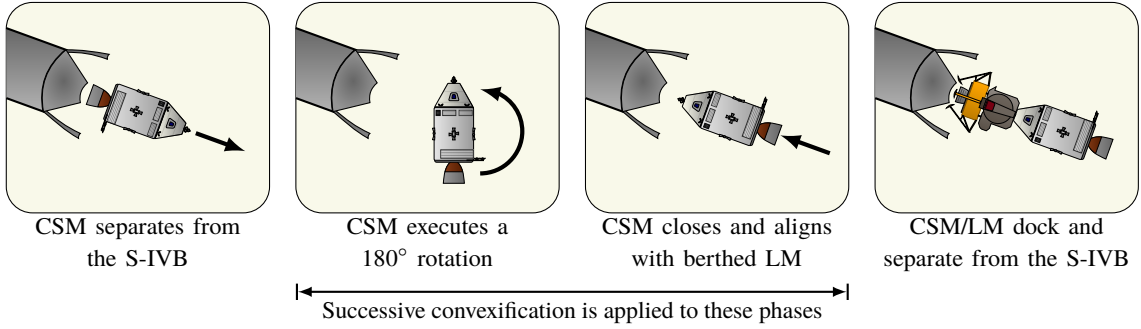

\begin{figure}
  \centering
    \includegraphics[width=\textwidth]{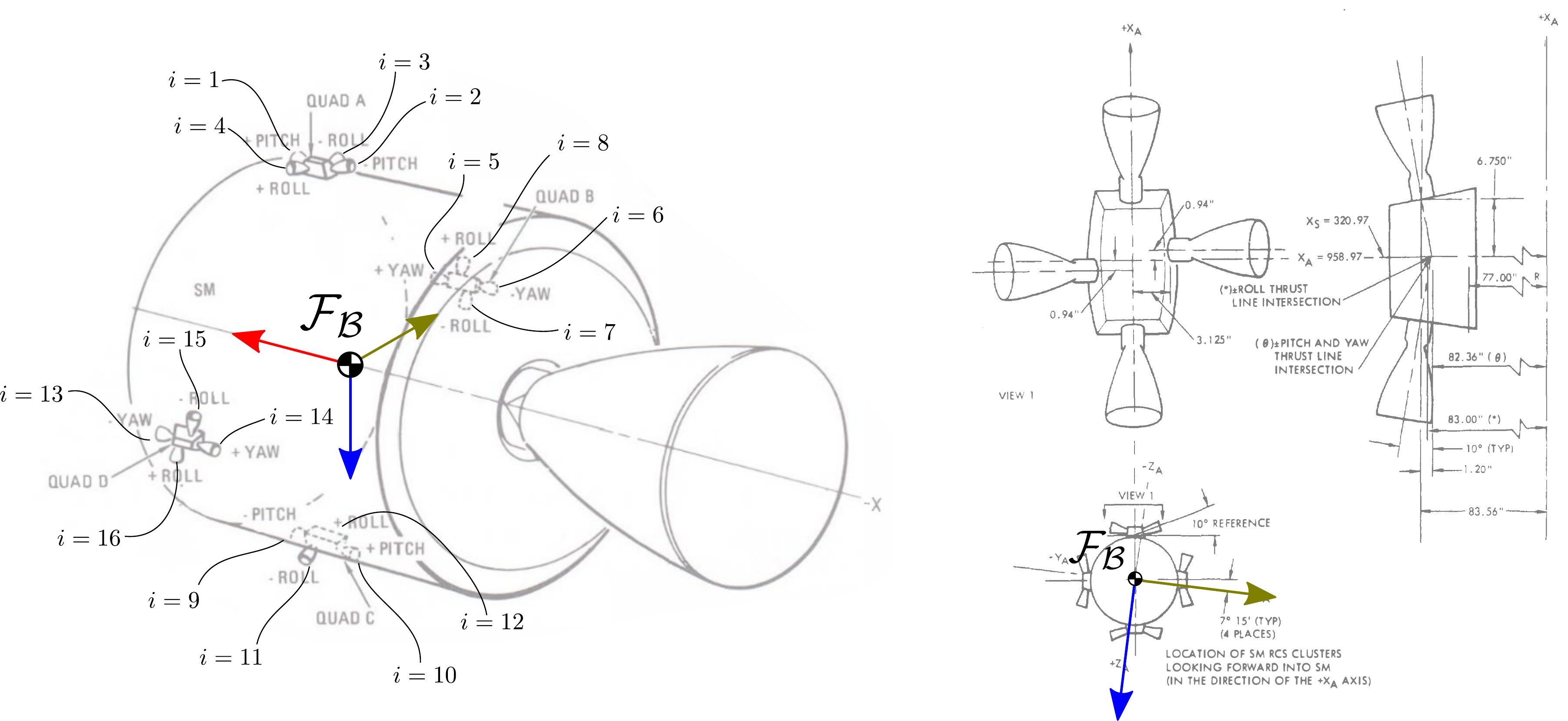}
  \caption{Layout of the Apollo SM RCS \cite[Figure~2.5-1]{aoh_vol1}
    \cite[Figure~2-5]{apollo_mass}, composed of four ``quads'' with four
    independent hypergolic pressure-fed thrusters. Our body frame and index
    number assignment to each thruster are overlaid. The figure on the right
    shows the 10$^\circ$ cant angle, the 7$^\circ$ 15$'$ angular offset and the
    4.8~cm offset in the quad and thruster layouts.}
  \label{fig:sm_rcs}
\end{figure}

We consider the Apollo CSM transposition and docking maneuver
\cite[Section~2.13.1.1]{aoh_vol1}. As illustrated in
Figure~\ref{fig:transposition_docking_illustration}, this maneuver uses the
\defabbrev{SM}{Service Module} RCS thrusters to dock with the LM, which is
housed inside the S-IVB third stage. The maneuver takes place after translunar
injection and we consider in particular the transposition and docking phases
(i.e. initial separation and final extraction are ignored, but can be handled
similarly).

The SM RCS system is composed of four similar, independent ``quads'' located
90$^\circ$ apart about the SM circumference, as illustrated in
Figure~\ref{fig:sm_rcs}. Each quad is composed of four independent hypergolic
pressure-fed pulse-modulated thrusters, yielding a total of $M=16$ control
inputs \cite[Section~2.5.1]{aoh_vol1}. Each quad is offset by $7^\circ~15'$ from
the spacecraft $y$ and $z$ body axes. Furthermore, the thrusters are canted
$10^\circ$ away from the SM outer skin to reduce the effects of exhaust gas on
the SM structure. Additionally, the two roll engines on each quad are
offset-mounted to accommodate plumbing. The document \cite{apollo_mass} provides
the complete geometry necessary to model the CSM and its RCS subsystem to high
fidelity. The CSM mass and inertia at the start of transposition and docking are
specified in \cite[Table~3.1-2]{apollo_mass}:
\begin{equation}
  m \approx 30322.9~\text{kg},\quad
  J \approx
  \begin{bmatrix}
    49248.7  &   2862.1 &   -370.1 \\
    2862.1 & 108514.2  &  -3075.0 \\
    -370.1 &  -3075.0  & 110771.7
  \end{bmatrix}~\text{kg}\cdot{m}^2.
\end{equation}

The expected RCS propellant consumption during transposition and docking for a
G-type mission (i.e. Apollo~11) was 50~kg, which is a less than 0.2~\% change in
mass. As a result, we argue that in this case ignoring mass depletion effects is
reasonable for trajectory planning.

The RCS thrusters are capable of producing approximately $\|\hat f_i\|_2=445$~N
of thrust in steady-state operation
\cite[Figure~2.5-8,~Table~4.3-1]{aoh_vol1}. However, during transposition and
docking they are pulse-fired in bursts \cite{apollo_news}
\cite[Section~2.5.1]{aoh_vol1}. In this mode of operation, the minimum electric
on-off pulse width is 12~ms \cite[Section~2.5.2.3.1]{aoh_vol1}, generating an
irregular burst of thrust that lasts for upwards of 65~ms and with a peak of 300
to 350~N \cite[Figure~2.5-9]{aoh_vol1}. Accordingly, we set
$\Delta t_{\min}=100$~ms and $\Delta t_{\max}=500$~ms, where the lower bound
somewhat avoids operating in the highly irregular thrust profile region of the
very short pulse widths. For the scope of this work, irregularity in the thrust
profile is assumed to be corrected for by a feedback control scheme that will
track our open-loop trajectory generated via successive convexification. We
note, however, that as long as a thrust profile can be expressed as a function
of time and pulse width, the irregularity can also be handled by successive
convexification.

Table~\ref{tab:parameters} summarizes the numerical values that were used to
obtain the results presented in the following section. We assume that the
transposition and docking maneuver is rotationally rest-to-rest with the CSM and
LM positive $x$ axes initially collinear, such that the CSM must perform a
180$^\circ$ flip as part of the maneuver. The docking impact velocity is set to
0.1~m/s, which is within the
specifications~\cite[Section~3.8.2.3]{csm_aoh}. Similarly, the propagation error
over $[0,t_f]$ must satisfy the required docking tolerances of 0.3~m in
position, 0.15~m/s in velocity, 10$^\circ$ in angular alignment and 1 $^\circ$/s
in angular velocity. Note that the docking geometry has the CSM assume an
approximately 60$^\circ$ roll with respect to the LM
\cite[Figure~2-4]{apollo_mass}. Instead of giving the desired terminal state of
the CSM, we give the LM state (position $p_l$, velocity $v_l$, attitude $q_l$,
and angular velocity $\omega_l$) in Table~\ref{tab:parameters}. This information
can be used to derive the desired CSM terminal state.

For the results in the following section, we have implemented Algorithm~\ref{alg:ftr}
in Python 3.7.2 using CVXPY 1.0.24 \cite{cvxpy} and MOSEK 9.0.87
\cite{mosek}. We run the algorithm on a Ubuntu~18.04.2 machine with an Intel
Core i7-6850K 3.6~GHz CPU and 64~GB RAM. Our implementation is publicly
available and the reader is encouraged to peruse the source code in complement
with this
paper\footnote{\url{https://github.com/dmalyuta/successive_rendezvous}.}.

\begin{table}
  \centering
  \begin{tabular}{r|l}
    Parameter & Value \\ \hline\hline
    $\Delta t_{\min}$ & 100~ms \\
    $\Delta t_{\max}$ & 500~ms \\
    $t_c$ & 2~s \\
    $\|\hat f_i\|_2$ & 445~N \\
    $r_a$ & 4~m \\
    $\Delta\theta_{\max}$ & 2$^\circ$ \\
    $\gamma$ & 30$^\circ$ \\
    $w_{tr}$ & $10^3$ \\
    $w_{vc}$ & $10^7$
  \end{tabular}
  \hspace{1cm}
  \begin{tabular}{r|l}
    Parameter & Value \\ \hline\hline
    $p_0, p_l$ & $(0,0,0)$~m, $(20,0,0)$~m \\
    $v_0, v_l$ & $(0,0,0)$~m/s, $(0,0,0)$~m/s \\
    $q_0, q_l$ & $(0,0,1,0)$, $(0,0,1,0)$ \\
    $\omega_0, \omega_l$ & $(0,0,0)$~deg/s, $(0,0,0)$~deg/s \\
    $p_{e,\max}$ & 1~cm \\
    $v_{e,\max}$ & 1~mm/s \\
    $\theta_{e,\max}$ & 0.5$^\circ$ \\
    $\omega_{e,\max}$ & 0.01~$^\circ$/s \\
    $\Delta J_{\text{tol}}$ & 0
  \end{tabular}
  \caption{Numerical parameters for Apollo CSM/LM transposition and docking.}
  \label{tab:parameters}
\end{table}

\subsection{Computed Rendezvous Trajectories}
\label{subsec:results}

We compute the rendezvous trajectory for durations of
$t_f\in\{150,250,350,450\}$~s. The total fuel consumption is computed as the
following integral:
\begin{equation}
  \label{eq:1}
  m_{\text{fuel}}(t_f) = m_{\text{fuel}}(0)+\int_0^{t_f} c_1 n_{\text{thruster}}(\tau)^2\dd\tau,
\end{equation}
where $n_{\text{thruster}}(\cdot):[0,t_f]\to\integer_{\ge 0}$ is the number of
thrusters firing at time $t$ and $c_1$~kg/s is the fuel consumption rate when a
single thruster is being fired \cite[Section~4.3.4.1.2]{csm_aoh}. As an
approximation, we use $c_1=0.168$~kg/s which is the fuel consumption in
steady-state operation \cite[Table~4.3-1]{csm_aoh}.

\begin{figure}
  \centering
  \includegraphics[width=0.6\textwidth]{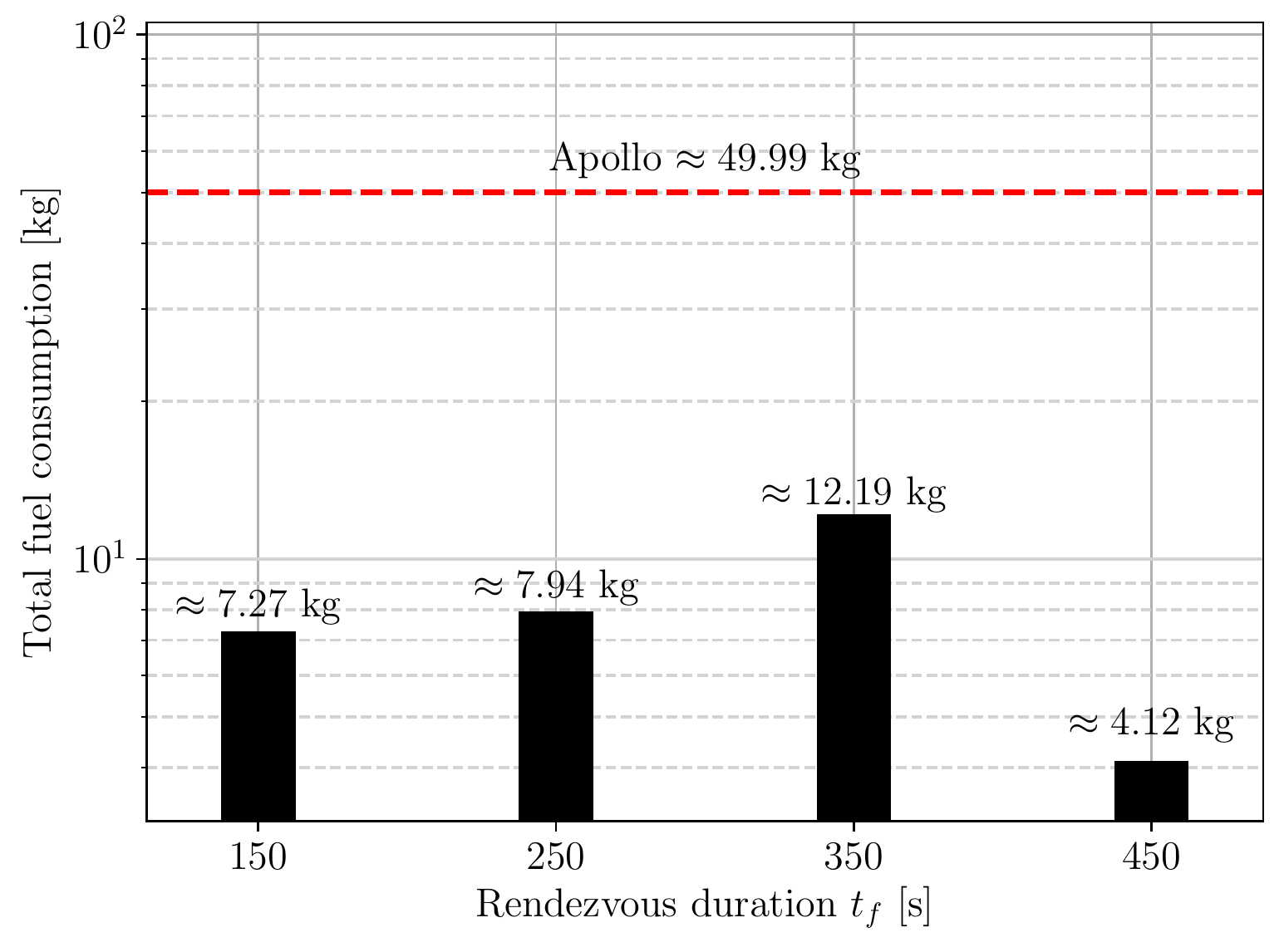}
  \caption{Total fuel consumption of the converged trajectories output by
    Algorithm~\ref{alg:ftr}. The Apollo G-type mission nominal fuel consumption
    is shown as reference.}
  \label{fig:cost_vs_tf}
\end{figure}

Figure~\ref{fig:cost_vs_tf} illustrates the resulting total fuel consumption of
the rendezvous trajectories output by Algorithm~\ref{alg:ftr}. First of all, we
note that the converged trajectories have vastly superior fuel consumption
compared to the Apollo G-type mission design target. Over $90~\%$ fuel may be
saved using the trajectory obtained for $t_f=450$~s. Note that it is
theoretically expected that, based on our dynamics \eqref{eq:vanilla_dynamics},
the total minimum fuel consumption should decrease as $t_f$ increases, since the
CSM can fire shorter pulses to achieve the same
goal. Figure~\ref{fig:cost_vs_tf} does not exhibit this trend because we are
finding only (roughly) locally optimal trajectories and not globally optimal
ones.

\begin{figure}
  \centering
  \includegraphics[width=0.9\textwidth]{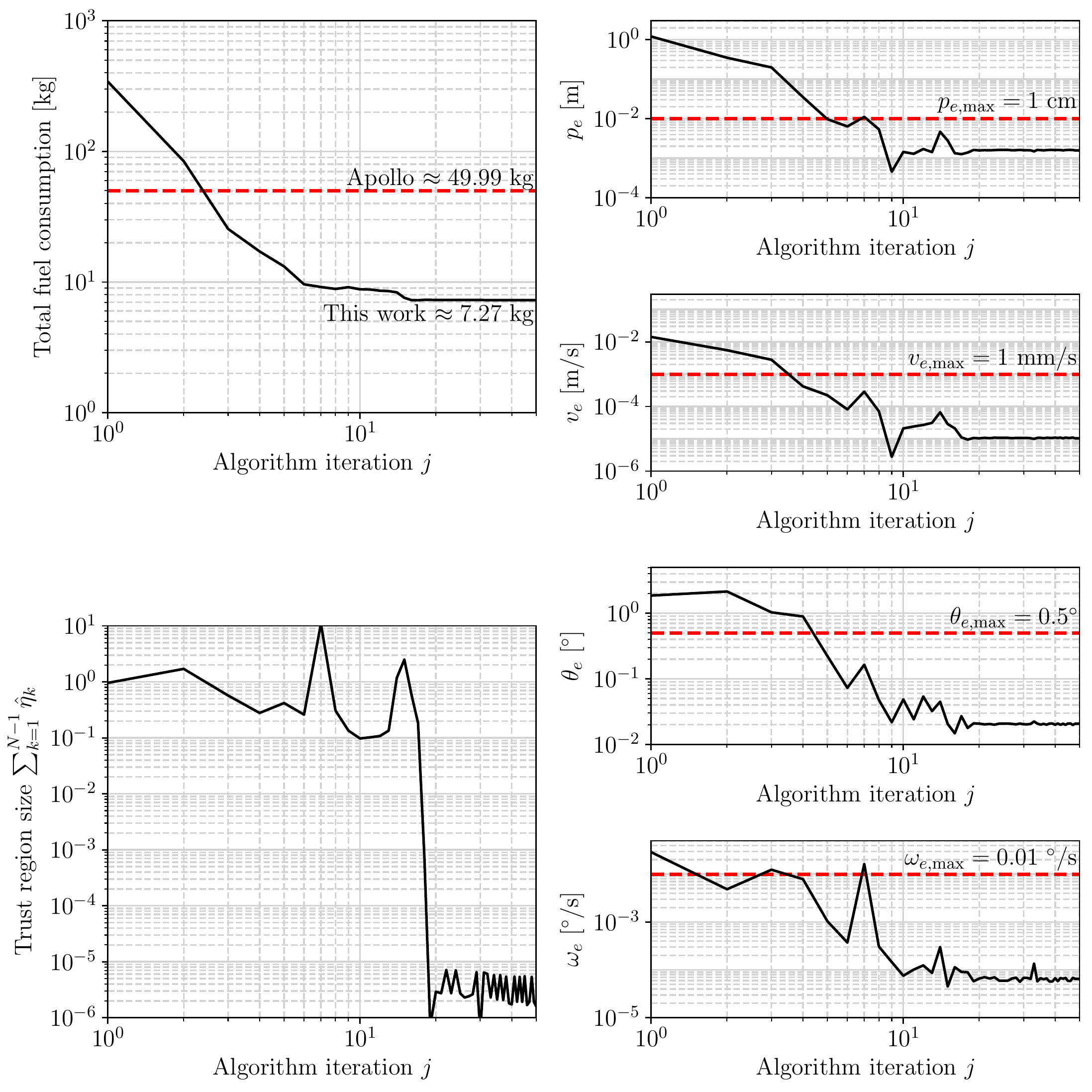}
  \caption{Algorithm~\ref{alg:ftr} convergence process. It takes only a few
    iterations to get a more fuel-optimal trajectory than the Apollo G-type
    mission design target. After that, it takes only a few more iterations
    before a feasible trajectory is found (i.e. the trust region decreases to a
    quasi-zero value and the propagation error over $[0,t_f]$ falls to within
    the required tolerance).}
  \label{fig:progress}
\end{figure}

Figure~\ref{fig:progress} provides the progress of Algorithm~\ref{alg:ftr}
across the successive convexification iterations for the case of $t_f=150$~s. We
can see that the fuel cost decreases quasi-monotonically. It takes only three
iterations before the total fuel consumed is less than the Apollo G-type mission
design target. After a few more iterations, the trust region falls sharply and
the propagation error decreases to within acceptable bounds.

\begin{figure}
  \centering
  \begin{subfigure}{1\textwidth}
    \centering
    \includegraphics[width=0.65\textwidth]{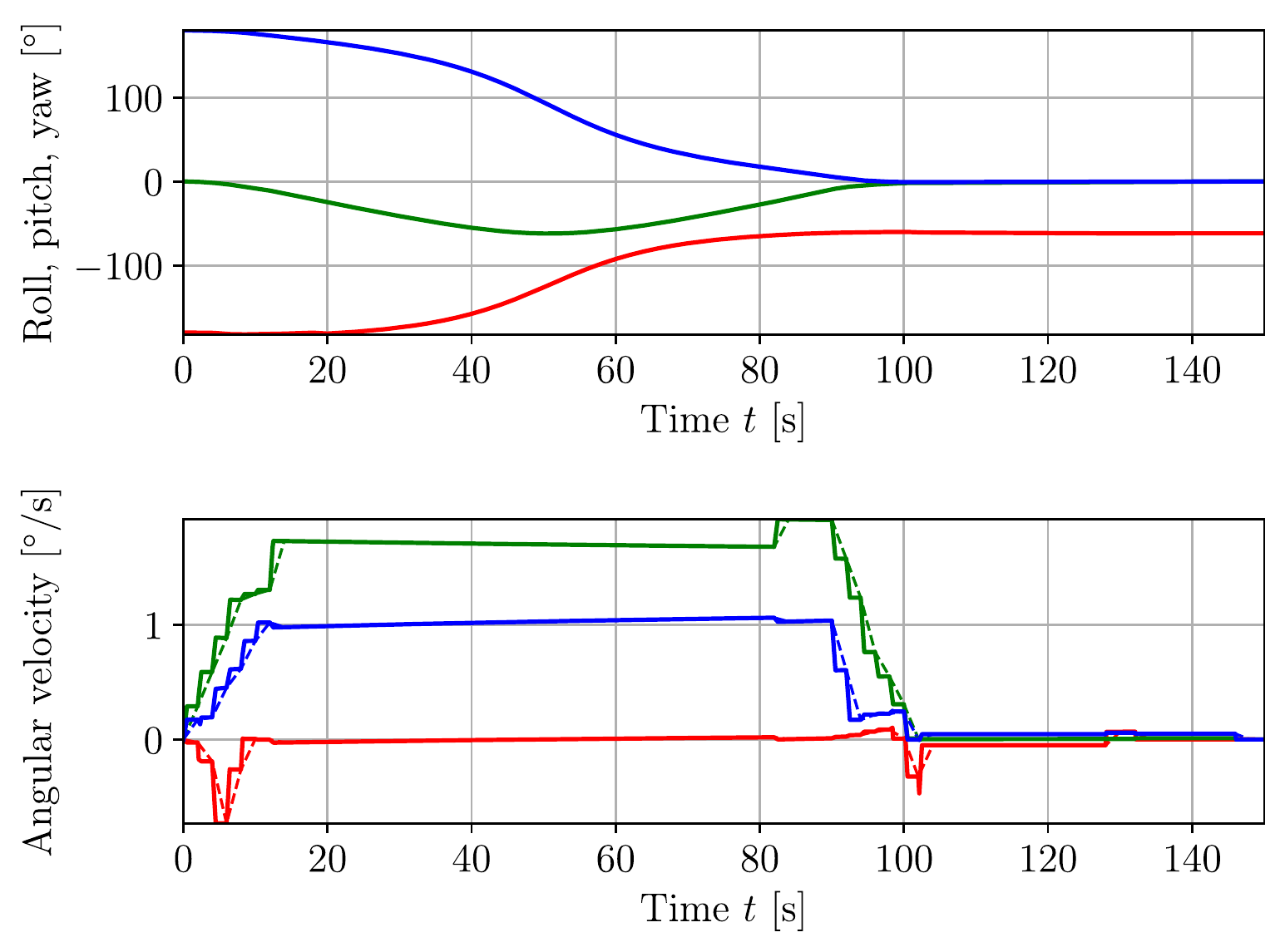}
    \caption{Converged CSM attitude trajectory.}
    \label{fig:attitude}
  \end{subfigure}
  \begin{subfigure}{1\textwidth}
    \centering
    \includegraphics[width=0.65\textwidth]{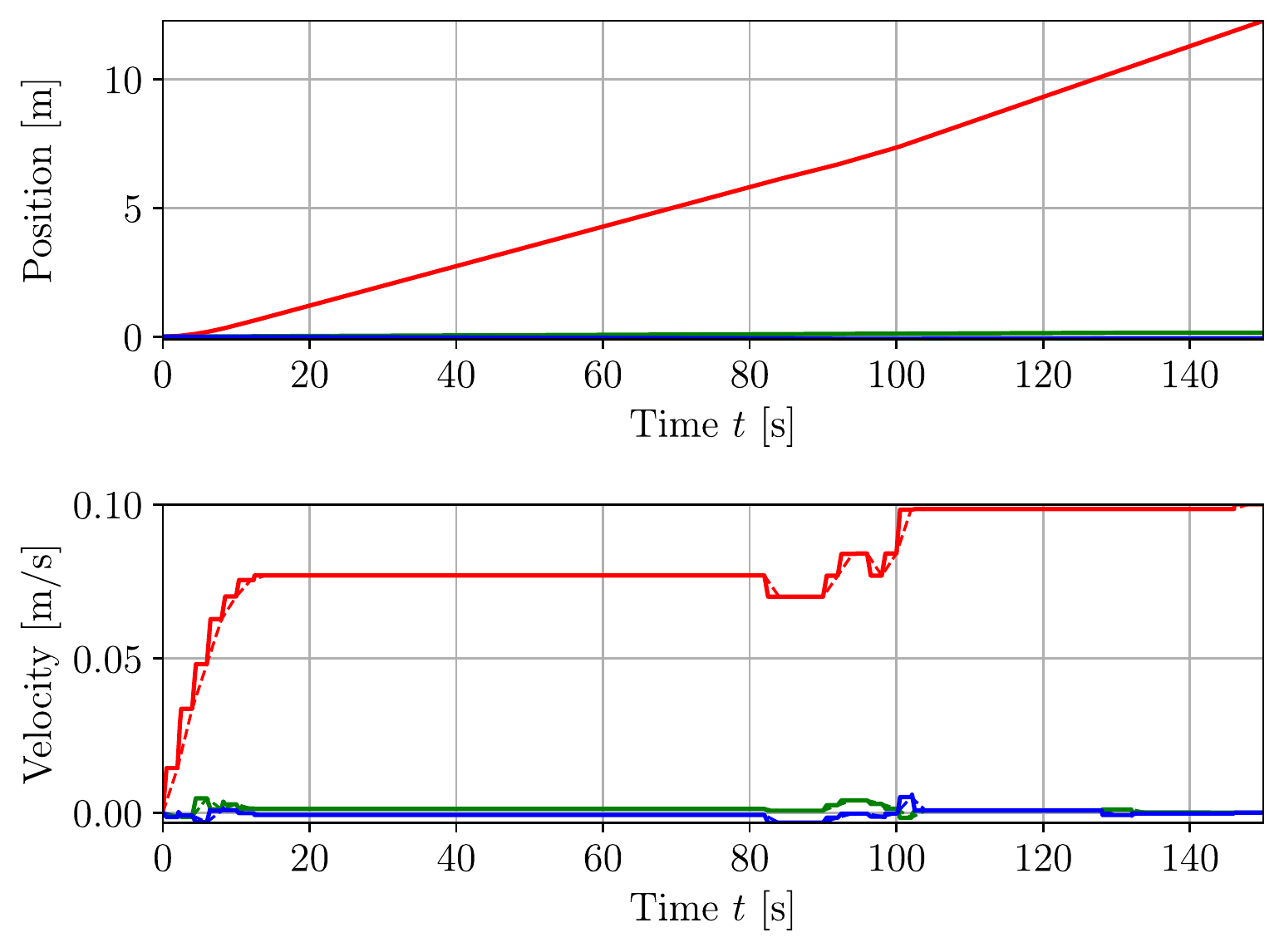}
    \caption{Converged CSM translation trajectory.}
    \label{fig:translation}
  \end{subfigure}
  \caption{Illustration of the converged trajectory obtained via
    Algorithm~\ref{alg:ftr} for a rendezvous duration of $t_f=150$~s. Roll,
    pitch and yaw correspond to the Tait-Bryan convention. The solid lines
    correspond to the continuous-time trajectory obtained by propagating the
    input trajectory through the non-linear
    dynamics~\eqref{eq:vanilla_dynamics}. The dashed lines correspond to the
    discrete-time trajectory output by the optimizer. Note the quasi-exact match
    of the two trajectories at the temporal grid nodes.}
  \label{fig:trajectory}
\end{figure}

Figure~\ref{fig:trajectory} shows the converged trajectory for $t_f=150$~s. The
total solver time, calculated as the cumulative time consumed by the core step
of solving Problem~\ref{problem:lrtop} (see
Figure~\ref{fig:successive_convexification_block_diagram} and
line~\ref{alg:ftr:line:solve} of Algorithm~\ref{alg:ftr}), is 6.8~seconds. Note
that, as required, the roll in Figure~\ref{fig:attitude} goes to a value of
$-60^\circ$, corresponding to the required CSM/LM relative docking
orientation. Also note that most of the maneuver, apart from translation along
the inertial $x$-axis, is finished after $100$~s, which corresponds to the time
when the CSM enters the $r_a=4$~m approach radius of the LM and where
constraints \eqref{eq:brtop_h} and \eqref{eq:brtop_i} require the forward
thrusters $i\in\{1,5,9,13\}$ to stay silent and the CSM to be within a
$2^\circ$ error angle of the terminal attitude.

\section{Future Work}
\label{sec:future}

Future work on this algorithm can proceed in several directions. First, we want
to better understand the state-triggered constraint locking property of
Section~\ref{subsec:stc} and to provide more rigorous remedies for it. Second,
the discretization operation of Section~\ref{subsec:propagation} is currently a
performance bottleneck because it is implemented in Python, unlike the C++
solver. This can be alleviated through a compiled language implementation and by
leveraging the fact that \eqref{eq:propagation_with_resetting} is parallelizable
across control intervals~\cite{Reynolds2019b}. Finally, it may be desirable to develop the algorithm to a point where a guaranteed real-time on-board implementation may be possible.

\section{Conclusion}
\label{sec:conclusion}

In this paper we have presented an algorithm, based on successive convexification, that is able to
generate fuel-optimal 6-DoF spacecraft rendezvous trajectories. Our method is
able to handle constraints that are important to rendezvous and docking operations using the state-triggered constraint paradigm, which is able to model discrete decisions within a continuous optimization framework. In particular, state-triggered constraints were used to model a
minimum RCS thruster pulse width, below which an RCS thruster is to remain
silent. State-triggered constraints were also used to model plume impingement constraints that were enforced only when the chaser vehicle lies within a specified radius of the target vehicle. Using the real world example of Apollo command and service module transposition and
docking with the lunar module, we showed that the method is able to quickly find trajectories which are up to $90\%$ more fuel optimal than the Apollo design targets. We believe that the algorithm may serve as a useful engineering tool for rapidly generating trajectories during mission design trade studies.

\section*{Acknowledgments}

This research was partially supported by the National Science Foundation
(CMMI-1613235) and by grant NNX17AH02A from the National Aeronautics and Space
Administration. Government sponsorship acknowledged.

\bibliography{references}

\end{document}